\documentclass[a4paper, 12pt]{article}
\usepackage{mathtools}
\usepackage{amssymb}
\usepackage{enumitem}
\usepackage{amsthm}
\usepackage{authblk}
\allowdisplaybreaks[1]
\usepackage{dsfont}
\usepackage{amssymb}
\usepackage{latexsym}
\usepackage{amsmath}
\usepackage{color}
\usepackage{comment}
\usepackage{amsthm}
\usepackage[pdftex]{graphicx}
\usepackage{layout} 
\usepackage{ulem}
\usepackage{bm}
\usepackage{bbm} 
\usepackage[bbgreekl]{mathbbol} 
\usepackage{authblk}
\usepackage{setspace} 
\usepackage{natbib}
\newif\ifrs
\rstrue
\ifrs \usepackage{mathrsfs} \fi  
\newif\ifcol
\coltrue 

\usepackage{etoolbox}

\newtheorem{theorem}{Theorem}[section]

\newtheorem{lemma}[theorem]{Lemma}

\newtheorem{proposition}[theorem]{Proposition}

\newtheorem{remark}[theorem]{Remark}

\numberwithin{equation}{section}
\newtheorem*{theorem*}{Theorem}
\newtheorem*{ass*}{Assumption}
\newtheorem*{lemma*}{Lemma}
\newtheorem*{definition*}{Definition}
\newtheorem*{proposition*}{Proposition}
\newtheorem*{corollary*}{Corollary}
\newtheorem*{remark*}{Remark}
\newtheorem*{example*}{Example}
\newif\ifcol
\colfalse
\ifcol
\newcommand{\colorr}{\color[rgb]{0.8,0,0}}

\newcommand{\colorn}{\color[rgb]{1,1,1}}

\else

\newcommand{\colorr}{\color{black}}

\newcommand{\colorn}{\color{black}}

\fi
\newif\ifcol
\coltrue
\ifcol

\else
\fi
\newif\ifcol
\colfalse
\ifcol

\else
\fi
\newif\ifcol
\colfalse
\ifcol

\else
\fi
\newif\ifcol
\colfalse
\ifcol

\else
\fi
\newif\ifcol
\colfalse
\ifcol

\else
\fi
\excludecomment{en-text}
\includecomment{jp-text}
\includecomment{comment}
\def\bd{\begin{description}}
\def\ed{\end{description}}

\def\D2{\bbD_{2,\infty-}}

\def\D{{\bf D}}

\def\cala{{\cal A}}
\def\calb{{\cal B}}
\def\calc{{\cal C}}

\def\cale{{\cal E}}
\def\calf{{\cal F}}
\def\calg{{\cal G}}
\def\calh{{\cal H}}

\def\call{{\cal L}}
\def\calm{{\cal M}}
\def\caln{{\cal N}}

\def\calp{{\cal P}}

\def\calv{{\cal V}}
\def\calw{{\cal W}}
\def\calx{{\cal X}}
\def\caly{{\cal Y}}

\def\ds{\displaystyle}
\def\yeq{\>=\>}
\def\yleq{\>\leq\>}

\def\sfd{{\sf d}}

\def\ep{\epsilon}
\def\half{\frac{1}{2}}

\def\cadlag{c\`adl\`ag\ }

\def\y{\vspace*{3mm}\\}
\def\halflineskip{\vspace*{3mm}}
\def\nn{\nonumber}
\def\be{\begin{equation}}
\def\ee{\end{equation}}
\def\bea{\begin{eqnarray}}
\def\eea{\end{eqnarray}}
\def\beas{\begin{eqnarray*}}
\def\eeas{\end{eqnarray*}}
\def\bi{\begin{itemize}}
\def\ei{\end{itemize}}
\def\im{\item}
\def\bd{\begin{description}}
\def\ed{\end{description}}
\newcommand{\bbA}{{\mathbb A}}
\newcommand{\bbB}{{\mathbb B}}

\newcommand{\bbD}{{\mathbb D}}

\newcommand{\bbF}{{\mathbb F}}
\newcommand{\bbG}{{\mathbb G}}

\newcommand{\bbN}{{\mathbb N}}

\newcommand{\bbR}{{\mathbb R}}

\newcommand{\bbT}{{\mathbb T}}

\newcommand{\bbV}{{\mathbb V}}
\newcommand{\bbW}{{\mathbb W}}

\newcommand{\ol}{\overline}

\newcommand{\wh}{\widehat}
\newcommand{\wt}{\widetilde}

\newcommand{\indep}{\perp\!\!\!\perp}
\newcommand{\oln}{\ol{n}}
\DeclareMathOperator*{\esssup}{ess\,sup}
\DeclareMathOperator*{\essinf}{ess\,inf}

\begin{document}

\title{Log-rank test with coarsened exact matching
\footnote{
This work was in part supported by 
Japan Science and Technology Agency CREST JPMJCR2115; 
the research project with Kitasato University; 
Japan Society for the Promotion of Science Grants-in-Aid for Scientific Research 
No. 
23H03354 (Scientific Research);  
and by a Cooperative Research Program of the Institute of Statistical Mathematics.
}
}
\author[1,2]{Tomoya Baba}

\author[1,2,3]{Nakahiro Yoshida}
\affil[1]{Graduate School of Mathematical Sciences, University of Tokyo
\footnote{Graduate School of Mathematical Sciences, University of Tokyo: 3-8-1 Komaba, Meguro-ku, Tokyo 153-8914, Japan. e-mail: nakahiro@ms.u-tokyo.ac.jp}
        }
\affil[2]{Japan Science and Technology Agency CREST
       }
\affil[3]{The Institute of Statistical Mathematics
       }
\maketitle
\noindent
\ \\
{\it Summary} \ 
It is of special importance in the clinical trial to compare survival times between the treatment group and the control group. Propensity score methods with a logistic regression model are often used to reduce the effects of confounders. However, the modeling of complex structures between the covariates, the treatment assignment and the survival time is difficult. In this paper, we consider coarsened exact matching (CEM), which does not need any parametric models, and we propose the weighted log-rank statistic based on CEM. We derive asymptotic properties of the weighted log-rank statistic, such as the weak convergence to a Gaussian process in Skorokhod space, in particular the asymptotic normality, under the null hypothesis and the consistency of the log-rank test. Simulation experiments are also conducted to compare the performance of the log-rank statistic with a propensity score method and CEM. Simulation studies show that the log-rank statistic based on CEM is more robust than the log-rank statistic based on the propensity score.

\section{Introduction}

In a clinical trial, %
suppose that each individual in a given collection $\bbG$ of individuals is assigned to one of groups $\bbG^z$ ($z\in\{0,1\}$) : 
the treatment group $\bbG^1$ and the control group $\bbG^0$. 
In this paper, we consider a comparison between the survival times of the individuals in $\bbG^1$ and $\bbG^0$, 
and assess the treatment effect through the survival data. 
The survival time of the individual $i$ is denoted by $T^i$. 
For a random variable $U^i \colon \Omega \to (0, \infty)$ denoting the censoring time for $i\in\bbG$, we set
\beas
\widetilde{T}^i = T^i \wedge U^i, \quad
\delta^i = 1_{\{ T^i \le U^i \}}
\eeas
for $i\in\bbG$. 
Moreover, it is supposed that each individual $i\in\bbG$ has a covariate $X^i$, that is a random variable taking values in $\bbR^\sfd$. 
Then, $\wt{T}^i$, $\delta^i$ and $X^i$ ($i\in\bbG$) are observable in the trial. 

Let $Z^i$ denote the assignment of the individual $i\in\bbG$, that is, $Z^i=z$ if $i\in\bbG^z$. 
The variables $Z^i$ are in general random. 
If the sequence $(Z^i)_{i\in\bbG}$ is generated by a Bernoulli trial, then the clinical trial is randomly controlled. 
In practice, however, quite possibly the distribution of the assignments may depend on the covariates the individuals have. 
In such a situation, it is necessary to adjust the distributions of the covariates between $\bbG^1$ and $\bbG^0$ 
to correctly compare the outcomes observed in the trial. 
Based on the observations 
$(X^i, \widetilde{T}^i, \delta^i)$, $i\in\bbG$, we test the hypothesis that there is no effect of the treatment. 
If the distributions of the covariates are the same for $\bbG^1$ and $\bbG^0$, 
one can use an ordinary log-rank test as a test statistic. 
However,  it does not make sense if the distributions are different, 
since the so assessed ``treatment effect'' may be biased by the difference between the distributions of covariates. 
It is often the case that patients of better conditions are apt to be selected than others in an invasive treatment, and 
this selection bias would affect the assessment of the treatment effect. 
If a perfectly randomized clinical trial is impossible ethically or for any other reason, then it is necessary to adjust the covariates' distribution 
to correctly measure the treatment effects. 

Propensity score methods, such as matching and inverse probability of treatment weighting (IPTW), 
are widely used to estimate treatment effects from observational data, including time-to-event data: see \cite{austin2014use}.
There is some research on causal inference in survival analysis. 
For example, \cite{austin2013performance} considers the Cox regression model and conducts a Monte Carlo simulation 
to compare the performance of propensity score methods to estimate the hazard ratio. 
\cite{xie2005adjusted} proposes a weighted Kaplan-Meier estimator and a log-rank test statistic with IPTW. 
The asymptotic theory of the Cox regression model with IPTW is developed by 
\cite{chen2001causal} and \cite{guo2021adjusted}. 
Propensity score methods need an estimation of the unknown propensity score, and logistic regression is commonly used. 
On the other hand, there are nonparametric matching methods. 
Although these methods have the advantage of not requiring any models for data, 
there is little mathematical research except for \cite{abadie2012martingale}.

In this paper, we propose a weighted log-rank test combined with a matching method,  
and study asymptotic properties of the test. 
%
%
%
This paper is finally aiming at proving the asymptotic behavior of the test statistic when $\#\bbG^0\gg\#\bbG^1$. 
Such a problem occurs in the studies of rare diseases and in advanced medicine, 
where the number of treated patients is much smaller than ones cured by ordinary treatments. 
For example, as mentioned below, the new therapy with heart sheets 
is the case.

This work was inspired by the joint research with Osaka University 
on the statistical analysis of 
the autologous stem cell-sheet transplantation therapy for treating cardiomyopathy, 
that was introduced by Professor Yoshiki Sawa and his group at the medical school of Osaka University. 
The authors appreciate Professor Shigeru Miyagawa (Osaka University, Graduate School of Medicine) and his research group members. 
%
Many thanks go to Professor Nobuhiko Okada (Kitasato University) with his group for their great supports for the joint research project, 
and also to 
Professor Masahiro Takeuchi (University of Tokyo), Dr.\>Takumi Suzuki (University of Tokyo) 
and Dr.\>Madoka Takeuchi (University of Tokyo) for their valuable comments and helps.

\section{Weighted log-rank statistic with CEM}
Coarsened exact matching (CEM) is a statistical method that classifies all individuals into matched and unmatched individuals: see \cite{iacus2011multivariate, iacus2019theory}.
The algorithm of CEM is as follows.
The groups $\bbG^z$ is given as $\bbG^z=\{i\in\bbG;\>Z_i=z\}$ for every $z\in\{1,0\}$. 
They are exclusive sets of individuals that are in general random. 
Set a partition $(a)_{a\in\bbA}$ of a measurable set 
$\calx$ of $\bbR^\sfd$. 
Suppose that $\bbA$ is a finite set and each element $a \in \bbA$ is measurable.
Input the covariates $X^i$ ($i\in\bbG$) 
and define $\calp$ by 
\beas 
\calp &=& \bigg\{(i,j)\in\bbG^1\times\bbG^0;\>(X^i,X^j)\in\bigsqcup_{a\in\bbA}(a\times a)\bigg\}.
\eeas
Let $G^z=\varpi_z\calp$ for $z\in\{1,0\}$, where $\varpi_z$ is the 
projection $\bbG^1\times\bbG^0\to\bbG^z$. 
For $i\in\bbG$, the processes $(N^i_t)_{t\in\bbR_+}$ and $(Y^i_t)_{t\in\bbR_+}$ are defined by 
\beas
N^i_t = 1_{\{ \widetilde{T}^i \le t, \>\delta^i = 1 \}}\quad\text{and}\quad
Y^i_t = 1_{\{ \widetilde{T}^i\ge t \}},
\eeas
respectively. 
Write $a^i$ for $a\in\bbA$ when $X^i\in a$.
We define the weight process $(w^i_t)_{t\in\bbR_+}$ for $i\in\bbG$ as follows: 
\bea\label{0410081303}
w^i_t 
&=& 
\left\{ 
\begin{array}{lc}
1 & \text{ if } i \in G^1\y
\left( \sum_{j\in\bbG^0} 1_{\{ X^j \in a^i \}} Y^{j}_t \right)^{-}\sum_{i'\in\bbG^1}1_{\{ X^{i'} \in a^i \}} Y^{i'}_t & \text{ if }i\in G^0\y
0 & \text{ if }i\not\in G^1\sqcup G^0,
\end{array} \right.
\eea
where $x^-=1_{\{x\not=0\}}x^{-1}$ for $x\in\bbR$. 
Remark that in (\ref{0410081303}), 
\beas
\left( \sum_{j\in\bbG^0} 1_{\{ X^j \in a^i \}} Y^{j}_t \right)^{-}\sum_{i'\in\bbG^1}1_{\{ X^{i'} \in a^i \}} Y^{i'}_t 
&=&
\left( \sum_{j\in G^0} 1_{\{ X^j \in a^i \}} Y^{j}_t \right)^{-}\sum_{i'\in G^1}1_{\{ X^{i'} \in a^i \}} Y^{i'}_t 
\eeas
if $i\in G^0$. 
%
%
Let
\begin{align*}
\mathcal{N}^i_t = \int_0^t w^i_s dN^i_s, \quad
\overline{\mathcal{N}}^z_t = \sum_{i\in G^z} \mathcal{N}^i_t, \quad
\mathcal{Y}^i_t = w^i_tY^i_t, \quad
\overline{\mathcal{Y}}^z_t = \sum_{i \in G^z} \mathcal{Y}^i_t.
\end{align*}
Remark that $(\caln^i_t)_{t\in\bbR_+}$ and $(\ol{\caln}^z_t)_{t\in\bbR_+}$ are no longer counting processes.

To assess the treatment effect, we consider the weighted log-rank statistic
\bea\label{0409241245}
\mathcal{W}_t 
&= &
\int_0^t K^n_s 
\{ (\overline{\mathcal{Y}}^1_s)^-d\overline{\mathcal{N}}^1_s
 -  (\overline{\mathcal{Y}}^0_s)^-d\overline{\mathcal{N}}^0_s \} \quad (t \in \mathbb{R}_+), 
\eea
where
\bea\label{0409241246}
K^n_s 
&=&  
\big\{(\overline{\mathcal{Y}}^1_0 + \overline{\mathcal{Y}}^0_0)
(\overline{\mathcal{Y}}^1_0\overline{\mathcal{Y}}^0_0)^-\big\}^{1/2}
(\overline{\mathcal{Y}}^1_s + \overline{\mathcal{Y}}^0_s)^-
\overline{\mathcal{Y}}^1_s \overline{\mathcal{Y}}^0_s W^n_s \quad (s \in \mathbb{R}_+)
\eea
and $W^n$ is a locally bounded predictable process with respect to the filtration $\bbF$ specified below.
When $W^n_s=1$, the weighted log-rank statistic $\calw_t$ becomes the usual log-rank statistic. 

\section{Main results}\label{section1}
%
Suppose that the random elements $T^i$, 
$X^i$ and $U^i$ are defined on a complete probability space $(\Omega,\calf,P)$. 
The completeness is not restrictive since completion of a probability space is always possible. 
Denote by ${\sf N}$ the collection of all $P$-null sets. 
To prove the main result, we need some results from the martingale theory.
We modify the arguments in \cite{fleming2005counting} to apply martingale theory 
to the survival analysis with covariates and treatment assignments.
We begin with computing the compensator of the process $N^i$.
Define the process $(\dot{N}^i_t)_{t\in\bbR_+}$ by
$\dot{N}^i_t  = 1_{\{ \widetilde{T}^i \le t,\> \delta^i = 0 \}}$ ($t \in \mathbb{R}_+$),  
and the filtration $(\mathcal{F}_t)_{t \in \mathbb{R}_+}$ 
by 
\beas
\calf_t 
&=& \sigma[X^i, Z^i,N^i_s, \dot{N}^i_s \mid i \in \bbG, s \le t]\vee\sigma[{\sf N}].
\eeas
Then the filtration $(\calf_t)_{t\in\bbR_+}$ is right-continuous: 
$\calf_t= \cap_{\epsilon > 0} \calf_{t+\epsilon}$. 
See \cite{bremaud1981point}. 

For $i\in\bbG$, denote by 
\beas 
\big(P^{T^i}(B\mid X^i=x, Z^i=z);\>B\in\bbB[\bbR_+],x\in\bbR^\sfd,z\in\{0,1\}\big)
\eeas
a regular conditional probability of $T^i$ given $(X^i,Z^i)=(x,z)$. 
See \cite{ikeda1989stochastic} for the definition and basic properties of a regular conditional probability.
Let 
\bea\label{0409241416}
\Lambda^i(t \mid X^i = x, Z^i=z) 
&=& 
\int_0^t\frac{ P^{T^i}(ds \mid X^i = x, Z^i=z) }{ P^{T^i}([s, \infty) \mid X^i = x, Z^i=z)} 
\eea
and 
\beas
\Lambda^i(t \mid X^i,Z^i) &=& \Lambda^i(t \mid X^i = x,Z^i=z)\mid_{x = X^i,\>z=Z^i}.
\eeas

In what follows, we assume the following condition. 
\bd
\im[{\bf[C]}] 
\bd
\im[(i)] 
The random elements $(T^i,U^i)$ ($i\in\bbG$) are $\calf_0$-conditionally independent. 
\im[(ii)] 
For each $i\in\bbG$, $T^i$ and $U^i$ are $\calf_0$-conditionally independent. 
\im[(iii)]
For each $i \in \bbG$, 
$T^i$ and $\sigma[X^j, Z^j \mid j \not= i]$ are $\sigma[X^i, Z^i]$-conditionally independent.
\im[(iv)]
In $i\in\bbG$, the conditional law $\call\{T^i\mid X^i=x,Z^i=z\}$ is identical.
\ed
\ed

\begin{table*}
\caption{The causal diagram of the random variables.}\label{table1}%
\beas 
\begin{array}{ccccccc}
 & &Z^i&&&\hspace{-5mm}U^i&\\
     &\nearrow&&\searrow&&\hspace{6mm}\searrow&\\
X^i &             &\longrightarrow    &                 & T^i&\longrightarrow&\wt{T}^i\\
     &&&&&&\\
\end{array}
\eeas
\end{table*}
%
%
%
Table \ref{table1} shows the dependency of the $i$-th random variables.
For random variables $X$ and $Y$, 
the arrow $X \to Y$ denotes that the value of $X$ affects the value of $Y$.

Under Condition $[C]$ (iv), 
we can choose a common version $\Lambda(t\mid x,z)$ of $\Lambda^i(t\mid X^i=x,Z^i=z)$ such that 
\beas 
\Lambda^i(t\mid X^i=x,Z^i=z) &=& \Lambda(t\mid x,z)\qquad(t\in\bbR_+,\>x\in\bbR^\sfd,\>z\in\{0,1\})
\eeas
for all $i\in\bbG$. 
We consider the following conditions. 
\bd
\im[{\bf[A1]}] 
There exists a $\sigma$-finite measure $\rho$ on $([0,\tau],\bbB([0,\tau]))$ and a constant $L$ 
such that $\rho$ has no atom, 
$\rho([0,\tau])<\infty$, 
$\Lambda (\cdot\cap[0,\tau] \mid x,z)$ is dominated by $\rho$, and 
\beas
\big|p(s\mid x,z) - p(s\mid y,z) \big|&\leq& L\>|x - y|
\eeas
for all $x, y \in \calx$ and $z\in\{0, 1\}$, 
for a version of the Radon-Nikodym derivative $p(\cdot\mid x,z)=d\Lambda (\cdot \cap[0,\tau] \mid x,z)/d\rho$. 
Moreover, 
\bea\label{0410130535}
\|p\|_\infty \coloneqq \sup_{s\in[0,\tau],x\in\calx,z\in\{0,1\}}p(s\mid x,z)<\infty.
\eea 
\ed
Usually, $\rho$ is the Lebesgue measure on $[0,\tau]$ or some measure equivalent to it. However, logically we may consider a measure like the Cantor distribution. 
The continuity of the $\Lambda ([0,\cdot], \mid x,z)$ can be relaxed but this condition is adopted for simplicity of the analysis.

We regard a parameter $n\in\bbN$ 
as the driver of the asymptotic theory as $n\to\infty$. 
In what follows, $\bbG^z$, $G^z$ ($z\in\{0,1\}$) and $\bbA$ for a partition of $\calx$ will depend on $n$. 
Let $n_z=\# G^z$ for $z\in\{0,1\}$. 
The number $n_z$ is generally random and satisfies $n_z\leq\#\bbG^z$. 
The number $n$ is not necessarily equal to $\oln=n_0+n_1$. 
For example, $n=\#\bbG$ individuals enter the trial, 
and the number of individuals randomly or deterministically assigned to the treatment/control group after matching is $n_z$. 
The number $\ol{n}$ can be much less than $n$ due to matching. 
In this example, It is also possible to assume that the parameter $n$ only specifies $n_z$ as a function $n_z(n)$ of $n$, without corresponding to any number in reality.  

As suggested in Introduction, we will consider the situation where $\#\bbG^0\gg\#\bbG^1$. 
For example, it is the case where the number of observations is limited for the new therapy whose effect we want to know, while one can use big registry data obtained in the traditional treatments.  
In such a situation, the conditions become asymmetric between the two groups, as set as follows. 
\bd
\im[{\bf [A2]}] There exists a (bounded) measurable function $q=(q_s)_{s\in[0,\tau]}$ on $[0,\tau]$ such that 
\beas
\sup_{s\in[0,\tau]}\bigg|n_1^-\sum_{i\in G^1}Y^i_sp(s\mid X^i, Z^i)-q_s\bigg|
&\to^p& 
0
\eeas
as $n\to\infty$. 
\ed


\bd
\im[{\bf [A3]}] As $n\to\infty$, 
\beas
P\bigg[\inf_{i\in G^1}\sum_{j\in G^0}1_{\{X^j\in a^i\}}Y^j_\tau=0\bigg]&\to&0.
\eeas 
\ed

\bd
\im[{\bf [A4]}] $\sup_{t\in[0,\tau],i\in G^0}w^i_t\to^p0$ as $n\to\infty$. 
\ed

\bd
\item[{\bf[A5]}]
$n_1^- \to^p 0$ and 
$n_1^{1/2}d_n \to^p 0$ as $n\to\infty$, for $\displaystyle d_n = \max_{a \in \mathbb{A}} \text{diam}(a)$. 
\ed
%
%
\bd
\im[{\bf [A6]}] 
There exists a bounded measurable function $W^\infty=(W^\infty_s)_{s\in[0,\tau]}$ on $[0,\tau]$ such that 
\beas
\sup_{s \in [0, \tau]} |W^n_s - W^\infty_s| \to^p 0
\eeas
as $n\to\infty$. 
\ed

By constructing the log-rank statistic, we want to test the hypothesis
\bd
\im[{$\bf [H_0]$}] 
$\ds \Lambda(t \mid x,1) = \Lambda(t \mid x,0)$ $(t\in[0,\tau],\>x\in\calx)$. 
\ed

Before we consider asymptotic properties of the weighted log-rank statistic as $n \to \infty$, 
we should clarify the dependency of variables on $n$. 
The random elements $(X^i, Z^i, T^i, U^i)$ $(i \in \bbG)$ may depend on $n$.
The sets $\bbG$, $\bbG^z$ ($z = 0, 1$) and $\bbA$ depend on $n$ as we have mentioned before.
On the other hand, $\calx$, $\Lambda^i(t\mid X^i = x, Z^i = z)$ and $p(s\mid x, z)$ do not depend on $n$.

%
%
The following central limit theorem enables us to make a critical region of the test. 
The D-space on $[0,\tau]$ equipped with the Skorokhod topology is denoted by $D([0,\tau])$.
See \cite{jacod2003limit, billingsley1999convergence}. 
\begin{theorem}\label{0410120204}
Suppose that conditions $[C]$, $[A1]$-$[A6]$ and $[H_0]$ are fulfilled. Then 
$\calw\to^d\bbW$ in $D([0,\tau])$ as $n\to\infty$, 
where $\mathbb{W}=(\bbW_t)_{t\in[0,\tau]}$ is a centered Gaussian process with independent increments and 
${\rm Var} [\mathbb{W}_t] = \frac{1}{2}\int_0^t (W^\infty_s)^2q_s\,\rho(ds)$ for $t\in[0,\tau]$. 
\end{theorem}
%
The proof of Theorem \ref{0410120204} is in Section \ref{0410131051}. 
Given a consistent estimator $\calv_\tau$ of the asymptotic variance $\text{Var}[\bbW_\tau]$ of $\bbW_\tau$, 
based on $\calw$, the critical region of the test of $[H_0]$ is, for example, 
\beas 
\calc &=& 
\{\calv_\tau^{-1/2}\calw_\tau\geq z_\alpha\}, 
\eeas
where $z_\alpha$ is the upper $\alpha$-point of the standard normal distribution. 
The shape of the critical region is chosen depending on the alternative hypothesis. 
%
We can construct a consistent estimator $\calv_\tau$ as follows.

\begin{proposition}\label{prop4}
Let $\calv_\tau$ be the following estimator of $\text{Var}[\bbW_\tau]$ in Theorem \ref{0410120204}.
\[
\calv_\tau = \frac{1}{n_1}\sum_{i \in G^1} \frac{1}{2}\int_0^\tau (W^n_s)^2 \,dN^i_s.
\]
Under conditions $[C]$, $[A1]$, $[A2]$, $[A5]$ and $[A6]$, we have $\calv_\tau \to^p \text{Var}[\bbW_\tau]$.
\end{proposition}

The proof of Proposition \ref{prop4} is given in section \ref{section2}.

We consider the consistency of the log-rank statistic $\mathcal{V}_\tau^{-1/2}\mathcal{W}_\tau$ from the convergence $\mathcal{W}_\tau \to^d \mathbb{W}$ obtained in Theorem \ref{0410120204}, 
where $\mathcal{V}_\tau$ is a consistent estimator of the variance $\text{Var}[\mathbb{W}_\tau]$.
Suppose that $[A1]$ and $[A6]$ hold.
Define the alternative hypothesis $[H_1]$ by

\begin{description}
\item[{$\bf [H_1]$}]
The function $p(t \mid x, z)$ in $[A1]$ satisfies
\[
p(t \mid x, 1) \le p(t \mid x, 0) \quad (t \in [0, \tau], x \in \mathcal{X}).
\]
Moreover, there exists $B \in \mathcal{B}([0, \tau])$ such that $\rho(B) > 0$ and
\[
\sup_{t \in B, x \in \mathcal{X}} W_t^\infty \{p(t \mid x, 1) - p(t \mid x, 0)\} < 0
\]
for the process $W^\infty$ in $[A6]$.
\end{description}

We consider the following condition to show the consistency of the log-rank statistic.

\begin{description}
\item[{$\bf [A7]$}]
There exists a positive number $\widetilde{q}_\tau$ such that
\[
n_1^- \sum_{i \in G^1} Y^i_\tau \to^p \widetilde{q}_\tau
\]
as $n \to \infty$.
\end{description}

\begin{theorem}\label{0510271447}
Suppose that conditions {$[A1]$-$[A7]$} 
and $[H_1]$ hold.
Let $\calv_\tau$ be a consistent estimator of $\text{Var}[\bbW_\tau]$ in Theorem \ref{0410120204}.
Then, we have
\begin{equation}\label{eq15} 
P(\mathcal{V}_\tau^{-1/2} \mathcal{W}_\tau < -R) \to 1 \quad (n \to \infty)
\end{equation}
for all $R > 0$.
\end{theorem}
See Section \ref{0510271450} for the proof of Theorem \ref{0510271447}.

%
%
%
%

\section{Sufficient conditions for the assumptions of the main theorems}
When $(X^i, Z^i, T^i, U^i)$ ($i \in \mathbb{G}$) is an i.i.d.\@ sequence, conditions $[A1]$-$[A7]$ can be simplified.
Let $\mathbb{G} = \{1, \ldots, n\}$.
Suppose that the sequence of random elements $\{(X^i, Z^i, T^i, U^i)\}_{i \in \bbG}$ is i.i.d.\@ for each $n \in \bbN$.
On the other hand, suppose that the conditional distribution $\mathscr{L}\{ (X^1, T^1, U^1) \mid Z^1 = z \}$ 
does not depend on $n$ for each $z \in \{0,1\}$.
Let $\calx \subset \bbR^d$ be compact.
For sequences $(a_n)$ and $(b_n)$, we write $a_n \lesssim b_n$ (resp.~$a_n \gtrsim b_n$) 
if there exists some constant $C > 0$ such that 
$a_n \le C b_n$ (resp.~$a_n \ge C b_n$) for all $n \in \mathbb{N}$.
Let $a_n \asymp b_n$ denote that $a_n \lesssim b_n$ and $b_n \lesssim a_n$ hold.
We consider the following conditions.
\begin{description}
\item[{\bf [B1]}] Constants $\beta, \theta \in (0, 1)$ and $\overline{d} \in \{ 1, \ldots, d \}$ 
satisfy $\beta < 2\theta$ and $\overline{d} \theta < 1$.

\item[{\bf [B2]}] 
$P(Z^1 = 1) \asymp n^{\beta - 1}$.

\item[{\bf [B3]}] The conditional distribution $\mathscr{L}\{ X^1 \mid Z^1 = z \}$ has a density $f(\cdot \mid z)$ 
with respect to a $\sigma$-finite measure $\nu$ on $\mathbb{R}^d$ with $\nu(\mathcal{X}) < \infty$ and it holds that 
\begin{align*}
&\esssup_{x \in \mathcal{X}} f(x\mid 1) < \infty, \\
&\essinf_{x \in \mathcal{X}} \left(f(x\mid 0)E[Y^1_\tau \mid X^1 = x, Z^1 = 0]\right) > 0, 
\end{align*}
where $\esssup_{x \in \mathcal{X}}$ and $\essinf_{x \in \mathcal{X}}$ are with respect to the measure $\nu$.

\item[{\bf [B4]}] $P(\widetilde{T}^1 = s, X^1 \in \calx \mid Z^1 = 1) = 0$ for all $s \in [0, \tau]$, 
and $P(\widetilde{T}^1 \ge \tau, X^1 \in \calx \mid Z^1 = 1) > 0$.

\item[{\bf [B5]}] $\displaystyle d_n = \max_{a \in \mathbb{A}} \text{diam}(a) = O(n^{-\theta})$ as $n \to \infty$ and 
$\displaystyle n^{-\overline{d}\theta} \lesssim \min_{a \in \mathbb{A}} \nu(a) 
\le \max_{a \in \mathbb{A}} \nu(a) \lesssim n^{-\overline{d}\theta}$.

\item[{\bf [A1$'$]}]
Condition $[A1]$ holds and the function $p(\cdot \mid x, 1)$ in $[A1]$ is continuous on $[0, \tau]$ for all $x \in \mathcal{X}$.
\end{description}

\begin{remark}\rm
\bd
\im[(i)]
We may regard $\overline{d}$ as the number of the continuous covariates. 
The situations we have in mind are as follows.
Suppose that the covariates $X^i$ ($i = 1, \ldots, n$) take values in 
$\mathbb{R}^{\overline{d}} \times \{0, 1\}^{d - \overline{d}}$. 
Let $m_L$ be the Lebesgue measure on $\mathbb{R}^{\overline{d}}$ 
and $m_C = \sum_{x \in \{0, 1\}^{d - \overline{d}}} \delta_x$, 
where $\delta_x$ is the Dirac measure on $\mathbb{R}^{d - \overline{d}}$. 
Define $\nu = m_L \times m_C$. 
Let $\mathcal{X} = [-1, 1]^{\overline{d}} \times \{0, 1\}^{d - \overline{d}}$ and 
$\mathbb{A}$ be the partition of $\calx$ defined by
\[
\mathbb{A} = \left\{ \mathcal{X} \cap \left(
\prod_{\ell = 1}^{\overline{d}} \left(\frac{k_\ell}{\lfloor n^\theta \rfloor},  \frac{k_\ell + 1}{\lfloor n^\theta \rfloor} \right] 
\times \prod_{m = \overline{d} + 1}^d \{ k_m \} \right) \;\middle|\; 
k_1 ,\ldots, k_d \in \mathbb{Z} \right\}.
\]
The partition $\mathbb{A}$ satisfies condition $[B5]$.
\im[(ii)] 
Under $[A1]$, in $[B4]$ the condition that 
$P(\widetilde{T}^1 = s \mid Z^1 = 1) = 0$ for all $s \in [0, \tau]$ can be replaced by 
the condition that 
$P(U^1 = s \mid Z^1 = 1) = 0$ for all $s \in [0, \tau]$. 

\ed
\end{remark}


For each $n \in \bbN$, let $V^i = (X^i, Z^i, T^i, U^i)$ $(i \in \bbG)$.
We have
\[
E[g(V^i) \mid X^i \in \calx, Z^i = 1] = \frac{ E[g(V^i) 1_{\{X^i \in \calx \}} \mid Z^i = 1] }{ P(X^i \in \calx \mid Z^i = 1)}
\]
for a nonnegative measurable function $g \colon \bbR^d \times \{0, 1\} \times (0, \infty) \times (0, \infty) \to \bbR$.
Note that $P(X^i \in \calx \mid Z^i = 1) > 0$ under $[B4]$.
We define 
\[
q_s = E[Y^1_s p(s \mid X^1, Z^1) \mid X^1 \in \calx, Z^1 = 1]\quad (s \in [0, \tau]).
\]
Consider a centered Gaussian process $\bbW = (\bbW_t)_{t \in [0, \tau]}$ with independent increments and 
$\text{Var}[\bbW_t] = \frac{1}{2} \int_0^t (W^\infty_s)^2 q_s \,\rho(ds)$.
We define the estimator $\calv_\tau$ of $\text{Var}[\bbW_\tau]$ as follows.
\[
\calv_\tau = \frac{1}{2} \int_0^t (W^\infty_s)^2\, dN^i_s.
\]

\begin{theorem}\label{0510271457}
Suppose that conditions $[C]$, $[A1']$, $[A6]$ and $[B1]$-$[B5]$ are satisfied.
Then, the following properties $(1)$-$(3)$ hold.
\begin{description}
\item[$(1)$] If $[H_0]$ is satisfied, then $\calw \to^d \bbW$ in $D([0, \tau])$ as $n \to \infty$.
\item[$(2)$] $\calv_\tau \to^p \text{Var}[\bbW_\tau]$ $(n \to \infty)$.
\item[$(3)$] If $[H_1]$ is satisfied, then $P\big( (\calv_\tau^{1/2})^- \calw_\tau < -R \big) \to 1$ as $n \to \infty$ for all $R > 0$.
\end{description}
\end{theorem}

See section \ref{0510271458} for the proof of Theorem \ref{0510271457}.

\section{Causal inference}
In the context of the causal inference, 
we consider a random field $\bbT^i:\Omega\times\{0,1\}\to(0,\infty)$ for $i\in\bbG=\bbG^1\sqcup\bbG^0$ on a sample space $\Omega$, 
to describe the potential outcome of the individual $i\in\bbG$, i.e., 
the random field $\bbT^i$ expresses counterfactual outcomes. 
When causal inference in question is survival analysis aiming at discovery of the treatment effect on survival time, 
$\bbT^i$ denotes the potential survival time of the individual $i$, 
and $T^i=\bbT(Z^i)$. 
When $i\in\bbG^z$, $\bbT^i(z)$ is observed, however, 
the counterfactual value of $\bbT^i(1-z)$ is never observed from the trial. 
In the causal inferential framework, one is interested in comparing, in some measure, 
the two conditional laws $\call\{\bbT^i(z)\mid X^i=x,Z^i=z\}$ and $\call\{\bbT(1-z)\mid X^i=x,Z^i=z\}$. 
Though $\bbT^i(1-z)$ is counterfactual when $Z^i=z$, $\call\{\bbT(1-z)\mid X^i=x,Z^i=z\}$ is estimable if the condition that 
\bd
\im[{\bf [N]}] 
$\call\{\bbT^i\mid X^i=x,Z^i=z\}=\call\{\bbT^i\mid X^i=x\}$ for $x\in\calx$ and $z\in\{0,1\}$. 
\ed
is met. 
Indeed, under $[N]$, $\call\{\bbT^i(1-z)\mid X^i=x,Z^i=z\}=\call\{\bbT^i(1-z)\mid X^i=x,Z^i=1-z\}
=\call\{T^i\mid X^i=x,Z^i=1-z\}$. 
Therefore, under the setting $[N]$, one can consider a test about equality between the effects $\bbT^i(1)$ and $\bbT^i(0)$ given $x$. 
Table \ref{table2} shows the relationship between variables and counterfactual outcomes.
It should be remarked that this test is of different nature from a test like $[H_0]$. 
To clarify the idea, consider a simple example such that $\bbT^i(z)=z+|X^i|+\delta^i$, where $\delta^i$ is a random variable having dependency on $Z^i$. 
Condition $[N]$ does not hold for this $\bbT^i$, so the comparison between $\bbT^i(1)$ and $\bbT^i(0)$ in the above sense is impossible. 
However, like the test of $[H_0]$, the comparison between $\call\{\bbT^i(1)\mid X^i=x,Z^i=1\}$ and $\call\{\bbT^i(0)\mid X^i=x,Z^i=0\}$ is still possible.
This is a comparison of distributions that is taking the effects of the unobserved confounder $\delta^i$ into account. 
In other words, if the property $[N]$ holds, then our test gives a causal inferential test of equality between 
$\call\{\bbT^i(1)\mid X^i=x,Z^i=z\}$ and $\call\{\bbT^i(0)\mid X^i=x,Z^i=z\}$ in some measure. 
\halflineskip
\begin{table*}
\caption{The causal diagram of the random variables with potential outcomes.}
\label{table2}
\beas 
\begin{array}{ccccccc}
 & &Z^i&&&\hspace{-5mm}U^i&\\
     &\nearrow&&\searrow&&\hspace{6mm}\searrow&\\
X^i &             &    &                 & T^i=\bbT^i(Z^i)&\longrightarrow&\wt{T}^i\\
     &\searrow&&\nearrow&&&\\
     &&\bbT^i&&&&
\end{array}
\eeas
\end{table*}

\section{Simulation Studies}
We conduct simulations to compare the performance of CEM and IPTW used in the log-rank test
and to confirm the consistency of the log-rank statistic with CEM. 
The first experiment aims to observe the robustness of CEM. 
We iterate data generation under the null hypothesis 300 times and 
draw histograms of the log-rank statistics with CEM and IPTW. 
While IPTW needs a parametric model to estimate the propensity score, 
CEM does not suppose any models when calculating the log-rank statistic. 
Thus, we consider two cases, with and without model misspecification. 
We generate two treatment assignments by using two kinds of logistic regression models, 
and we estimate the propensity score with only one of the logistic regression models.
In the second experiment, we iterate to generate data under the alternative hypothesis and 
calculate the log-rank statistics with CEM 300 times. 
We draw three histograms of the log-rank statistics based on data of different sample sizes to observe the consistency.

The detailed definition of the simulation procedures is as follows. 
We generate i.i.d.~samples $\{ (X^i, Z^i, \widetilde{T}^i, U^i) \}_{i = 1}^n$, 
where $n = 5000$ in the first experiment and $n = 2500$, $5000$ and $7500$ in the second experiment.
Let $X^i = (X^i_1, \ldots, X^i_5)$ ($i = 1, \ldots, n$) be $\mathbb{R}^3 \times \{0, 1\}^2$-valued random variables.
Suppose that $X^i_1, \ldots, X^i_5$ are independent for every $i = 1, \ldots, n$.
Let the first three variables $X^i_1$, $X^i_2$ and $X^i_3$ follow the standard normal distribution $N(0, 1)$ and 
$X^i_4$ and $X^i_5$ satisfy $P(X^i_2 = 1) = P(X^i_2 = 0) = 1/2$ for each $i = 1, \ldots, n$. 
We consider the following two models for generating $Z^i$.
\begin{description}
\item[Model (1)]
a logistic regression model without interaction terms
\[
P(Z^i = 1 \mid X^i) = \frac{\exp(-3.4 -0.2 (X^i_1 + \cdots + X^i_5))}
{1 + \exp(-3.4 -0.2 (X^i_1 + \cdots + X^i_5))}.
\]
\item[Model (2)]
a logistic regression model with interaction terms
\[
P(Z^i = 1 \mid X^i) = \frac{\exp(-3.7 -0.2 (X^i_1 + \cdots + X^i_5) + 0.5(X^i_1X^i_2 + X^i_1X^i_3))}
{1 + \exp(-3.7 -0.2 (X^i_1 + \cdots + X^i_5) + 0.5(X^i_1X^i_2 + X^i_1X^i_3))}.
\]
\end{description}
The mean of $\#\mathbb{G}^1$ for 300 data rounding to the nearest whole number is 
$143$ in Model (1) and $139$ in Model (2).
Since $\log 140 / \log 5000$ is approximately $0.58$, we may suppose that $\beta = 0.58$, 
where $\beta$ is the constant in $[B1]$. 
We generate the potential outcomes $\mathbb{T}^i(z)$ ($z = 0, 1$, $i = 1, \ldots, n$) by using the Cox hazards model.
For $x = (x_1, \ldots, x_5) \in \bbR^3 \times \{0, 1\}^2$ and $z \in \{0, 1\}$, 
let $\Lambda(t \mid x, z)$ ($t \ge 0$) denote the cumulative hazard function 
of the conditional distribution of $\bbT^i(z)$ given $X^i = x$.
In other words, define
\[
\Lambda(t \mid x, z) = \int_0^t \frac{P^{\bbT^i(z)}(ds \mid X^i = x)}{P^{\bbT^i(z)}([s, \infty) \mid X^i = x)}.
\]
If the condition $[N]$ holds, this function equals the function $\Lambda(t \mid x, z)$ in Equation (\ref{0409241416}).
We generate $\{\bbT^i(z)\}_{i = 1}^n$ ($z = 0, 1$) that follow the Cox hazard model
\[
\Lambda(t \mid x, z) = \begin{cases}
\displaystyle \int_0^t h_0(s) \exp(0.25(x_1 + \cdots + x_5))\,ds & (\text{under the null hypothesis})\\
\displaystyle \int_0^t h_0(s) \exp(- 0.4z + 0.25(x_1 + \cdots + x_5))\,ds & (\text{under the alternative hypothesis})
\end{cases}
\]
where $h_0(s)$ is the constant function $h_0(s) = e^{-2}$.
Let $T^i = Z^i\mathbb{T}^i(1) + (1 - Z^i)\mathbb{T}^i(0)$.
Suppose that $U^i$ ($i = 1, \ldots, n$) follows the uniform distribution $U(0, 10)$.

Let $\tau = 10$.
From the observations $\{ (X^i, Z^i, \widetilde{T}^i, \delta^i) \}_{i = 1}^n$, we calculate the log-rank statistics $\calw_\tau$ 
and the estimator of their asymptotic variances $\calv_\tau$ using CEM and IPTW.
Whether using CEM and IPTW, the definition of $\calw_\tau$ is given in (\ref{0409241245}) and (\ref{0409241246}).
Let $W^n = 1$.
The differences between the definition of the log-rank statistics with CEM and IPTW are 
the weight processes $w^i_s$ and the estimator of the asymptotic variance $\calv_\tau$.
First, we consider the case with CEM proposed in this paper.
Let $\mathcal{X} = [-5, 5]^3 \times \{0, 1\}^2$.
We define the partition $\mathbb{A}$ for CEM as follows.
\[
\mathbb{A} = \left\{ \mathcal{X} \cap \left(
\prod_{i = 1}^3 \left( \frac{5k_i}{\lfloor n^{0.3} \rfloor} ,\frac{5(k_i+1)}{\lfloor n^{0.3} \rfloor} \right] \times 
\{k_4\} \times \{k_5\} \right)
\;\middle|\; k_1, \ldots, k_5 \in \mathbb{Z} \right\}.
\]
Note that $\beta = 0.58 < 2\theta = 0.6 < 2/\overline{d} = 2/3$.
The process $w^i_s$ for the log-rank statistic with CEM is defined as (\ref{0410081303}), 
and the estimator of the asymptotic variance $\calv_\tau$ is given in Proposition \ref{prop4}.
Next, we define $w^i_s$ and $\calv_\tau$ for the log-rank statistic with IPTW.
This statistic is introduced in \cite{xie2005adjusted}.
We calculate $\widehat{w}^i$ as follows and let $w^i_s = \widehat{w}^i$ for all $s \ge 0$.
\begin{align*}
\widehat{w}^i &= Z^i/\widehat{p}^i + (1 - Z^i)/(1- \widehat{p}^i), \\
\widehat{p}^i &= \frac{\exp(\widehat{\alpha}_0 + \widehat{\alpha}_1 X^i_1 + \widehat{\alpha}_2 X^i_2)}
{1 + \exp(\widehat{\alpha}_0 + \widehat{\alpha}_1 X^i_1 + \widehat{\alpha}_2 X^i_2)}, 
\end{align*}
where $(\widehat{\alpha}_0, \widehat{\alpha}_1, \widehat{\alpha}_2)$ is the maximum likelihood estimator 
of the parameter $(\alpha_0, \alpha_1, \alpha_2)$ in the logistic regression model 
\[
P(Z^i = 1 \mid X^i) = \frac{\exp(\alpha_0 + \alpha_1 X^i_1 + \alpha_2 X^i_2)}
{1 + \exp(\alpha_0 + \alpha_1 X^i_1 + \alpha_2 X^i_2)}.
\]
Let 
\begin{align*}
\overline{N}_s &= \sum_{i \in \bbG} N^i_s, \quad
\overline{Y}_s = \sum_{i \in \bbG} Y^i_s,  \\
U_s &= \{\overline{\caly}^0_s (\overline{\caly}^1_s + \overline{\caly}^0_s)^-\}^2 \sum_{i \in \bbG^1}(w^i_s)^2Y^i_s +
\{\overline{\caly}^1_s (\overline{\caly}^1_s + \overline{\caly}^0_s)^-\}^2 \sum_{i \in \bbG^0}(w^i_s)^2Y^i_s
\end{align*}
for $s \ge 0$ and $i \in \bbG$
and we define
\begin{align*}
\calv_\tau = (\overline{\caly}^1_0 + \overline{\caly}^0_0)^- \overline{\caly}^1_0 \overline{\caly}^0_0
\int_0^\tau U_s \{\overline{Y}_s (\overline{Y}_s - 1)\}^- (\overline{Y}_s - \Delta\overline{N}_s)\,d\overline{N}_s.
\end{align*}

\begin{figure}[tbp]
\centering
\includegraphics[width =  119mm]{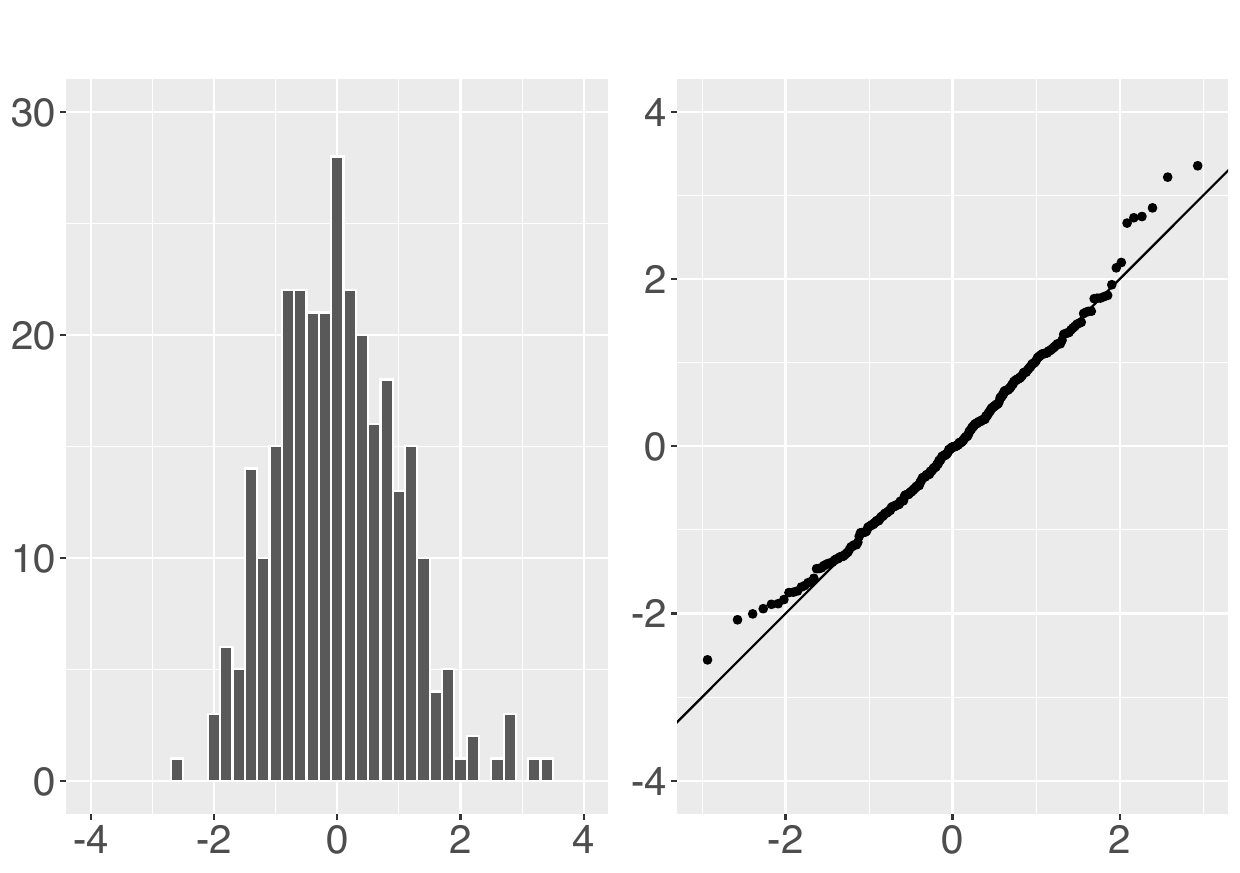}
\caption{Histogram and Q-Q plot of $\calv_\tau^{-1/2} \calw_\tau$ with IPTW 
under the null hypothesis and Model (1)}
\label{fig1}
\end{figure}
\begin{figure}[tbp]
\centering
\includegraphics[width = 119mm]{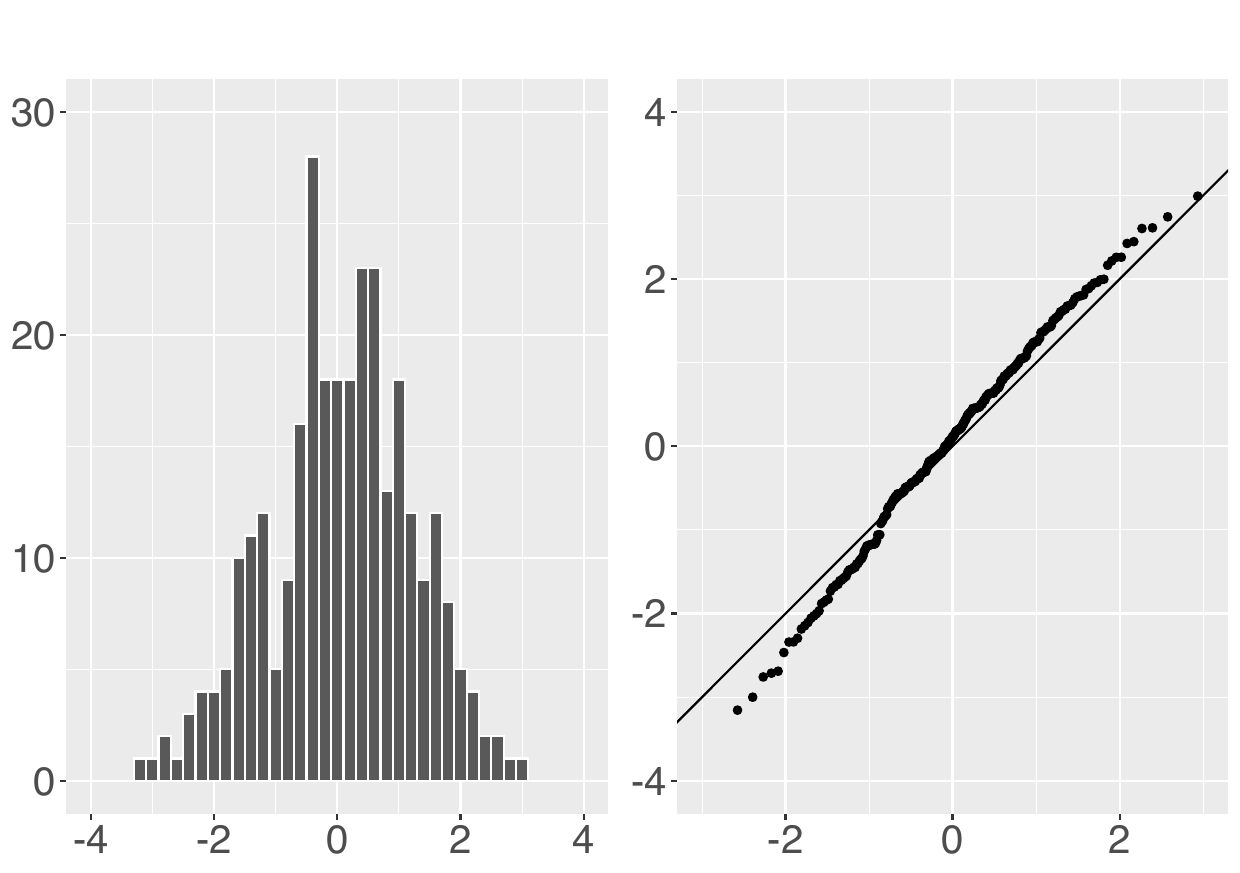}
\caption{Histogram and Q-Q plot of $\calv_\tau^{-1/2} \calw_\tau$ with CEM 
under the null hypothesis and Model (1)}
\label{fig2}
\end{figure}
\begin{figure}[tbp]
\centering
\includegraphics[width = 119mm]{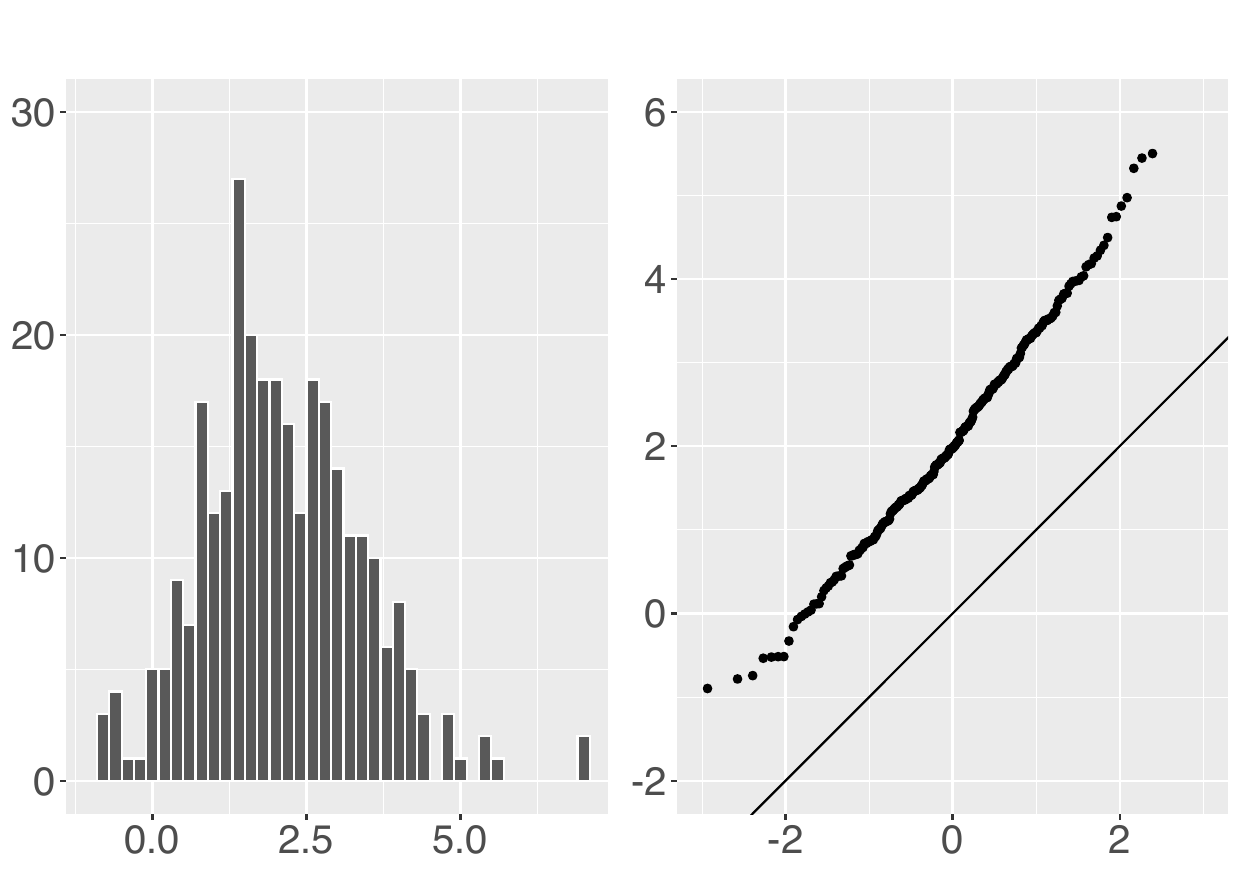}
\caption{Histogram and Q-Q plot of $\calv_\tau^{-1/2} \calw_\tau$ with IPTW 
under the null hypothesis and Model (2)}
\label{fig3}
\end{figure}
\begin{figure}[tbp]
\centering
\includegraphics[width = 119mm]{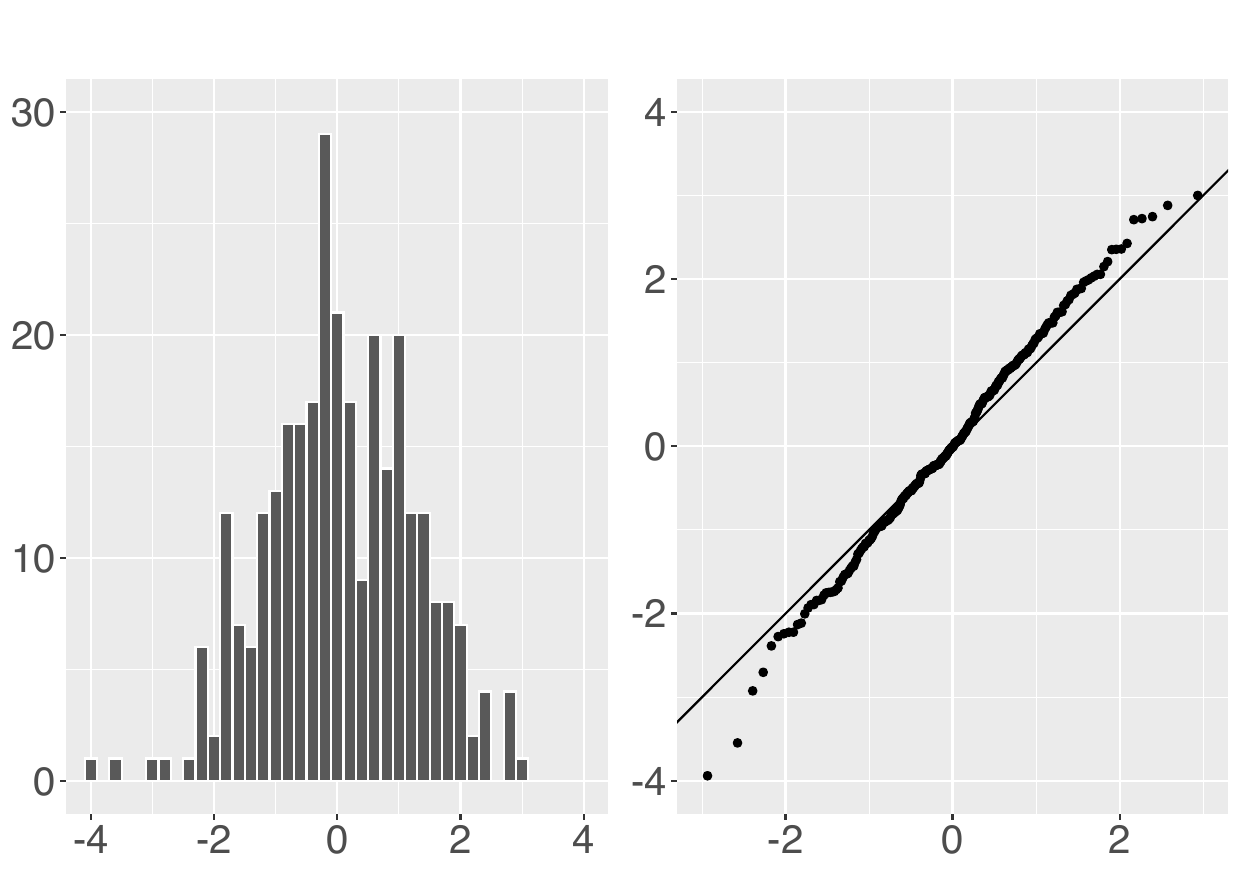}
\caption{Histogram and Q-Q plot of $\calv_\tau^{-1/2} \calw_\tau$ with CEM 
under the null hypothesis and Model (2)}
\label{fig4}
\end{figure}
Figures \ref{fig1}-\ref{fig4} show the result of the first experiment.
There is one outlier in -5.2 in \ref{fig2}.
From Figure \ref{fig1} and Figure \ref{fig2}, 
we can see that the performance of the log-rank statistics with CEM and IPTW is almost the same 
if we generate data that follows Model (1).
Moreover, the distribution of $\calv_\tau^{-1/2} \calw_\tau$ converges to the standard normal distribution.
On the other hand, Figure \ref{fig3} and Figure \ref{fig4} show that only the log-rank statistic with CEM works 
if we generate data that follows Model (2).
Therefore, we can conclude that the log-rank test with CEM performs well compared to the log-rank test with IPTW 
if the sample size of the control group is large enough and 
the relationship between the covariates and the treatment assignment is unknown.

\begin{figure}[tbp]
\centering
\includegraphics[width = 119mm]{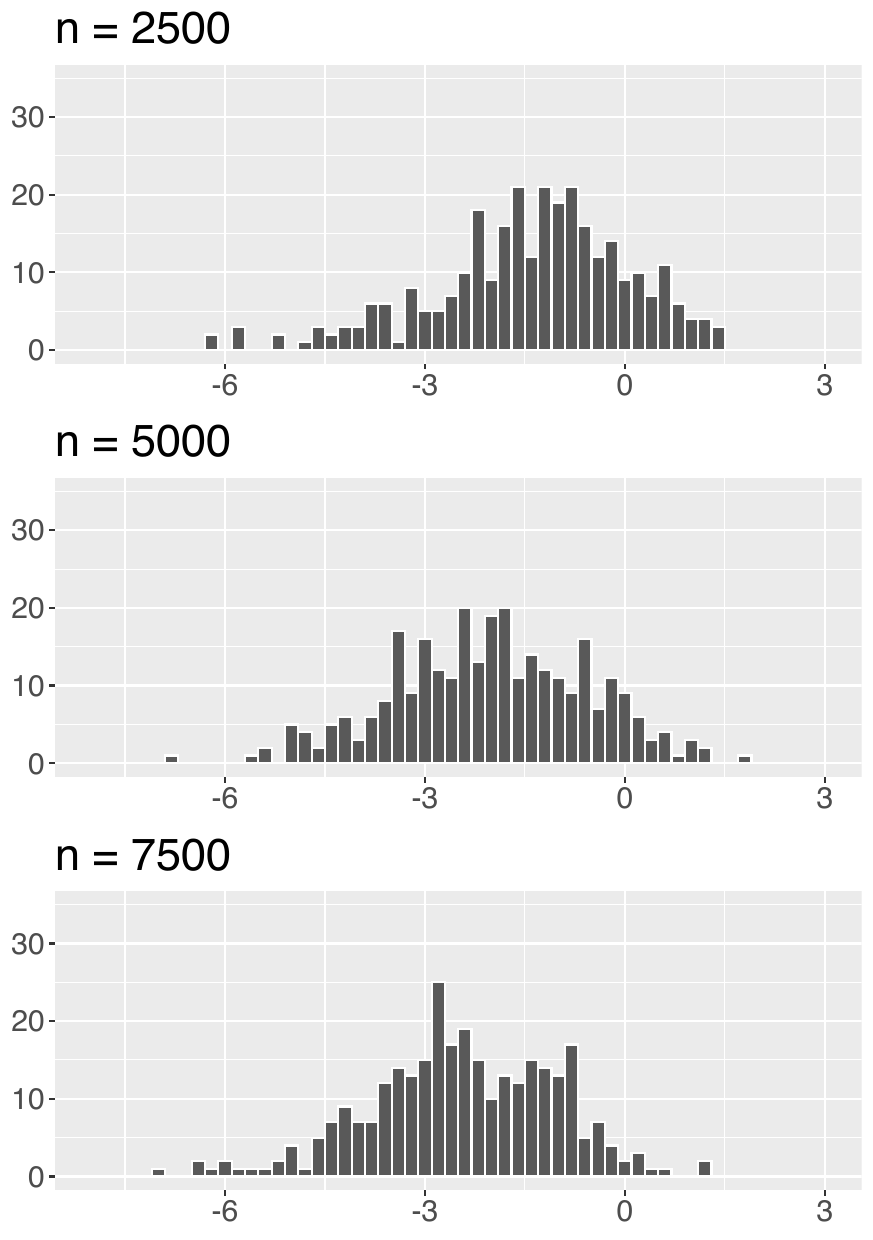}
\caption{Histograms of $\calv_\tau^{-1/2} \calw_\tau$ with CEM 
under the alternative hypothesis and Model (2)}
\label{fig5}
\end{figure}
Figure \ref{fig5} shows the result of the second experiment.
As the sample size increases, the histogram moves in the minus direction.
This result demonstrates the consistency of the log-rank test statistic with CEM.

\section{Proof of Theorem \ref{0410120204}}\label{0410131051}
For $i\in\bbG$, define the process $A^i=(A^i_t)_{t\in[0,\tau]}$ by 
\beas
A^i_t &=& \int_0^t Y^i_s \Lambda^i(ds \mid X^i, Z^i).
\eeas
%
%
\begin{proposition}\label{prop1} 
Suppose that $[C]$ is satisfied. Then 
$A^i$ is a (unique) compensator of $N^i$ with respect to $(\calf_t)_{t\in \mathbb{R}_+}$. 
\end{proposition}

Let $M^i_t=N^i_t-A^i_t$ for $t\in \mathbb{R}_+$ and $i\in\bbG$.

\begin{proposition}\label{2209250620} 
Suppose that $[C]$ is satisfied. Then 
$M^i$ is a square-integrable $(\calf_t)_{t\in \mathbb{R}_+}$-martingale for each $i\in\bbG$, and the predictable quadratic covariation process between $M^i$ and $M^{i'}$ is 
$\langle M^i,M^{i'}\rangle = 1_{\{i=i'\}}A^i$ for $i,i'\in\bbG$. 
\end{proposition}

See Appendix for the proof of Proposition \ref{prop1} and \ref{2209250620}.
Let
\beas
\mathcal{M}^z_t
&=& \sum_{i\in G^z}\int_0^t K^n_s(\overline{\mathcal{Y}}^z_s)^-w^i_s dM^i_s, 
\eeas
and 
\beas
\mathcal{E}^z_t
&=& 
\int_0^t K^n_s(\overline{\mathcal{Y}}^z_s)^-\sum_{i\in G^z}w^i_sY^i_s \Lambda(ds \mid X^i,z).
\eeas
Obviously, 
\beas
\mathcal{E}^z_t
&=& 
\int_0^t K^n_s(\overline{\mathcal{Y}}^z_s)^-\sum_{i\in G^z}\sum_{a\in\bbA}1_{\{X^i\in a\}}w^i_sY^i_s \Lambda(ds \mid X^i,z).
\eeas
Then, it holds that
\bea\label{0410121858}
\mathcal{W}_t 
&=& 
\mathcal{M}^1_t -\mathcal{M}^0_t + \mathcal{E}^1_t -  \mathcal{E}^0_t.
\eea



We will prove $\mathcal{M}^1$ converges in distribution to Gaussian process in $D[0, \tau]$ and 
$\mathcal{M}^0$ and $\mathcal{E}^1 - \mathcal{E}^0$ converges to 0 in probability 
under conditions $[C]$, $[A1]$-$[A6]$ and $[H_0]$.

For $z\in\{0,1\}$ and $i\in\bbG$, define the process $H^{z,i}=(H^{z,i}_s)_{s\in\bbR_+}$ by 
\bea\label{0410081534}
H^{z, i}_s = K^n_s (\overline{\caly}^z_s)^- w^i_sY^i_s\qquad(s\in\bbR_+)
\eea
for $z\in\{0,1\}$ and $i\in\bbG$. 
Define $\dot{\cale}^z_t$ and $\ddot{\cale}^z_t$ by 
\beas
\dot{\mathcal{E}}^z_t
&=&
\int_0^t 
\sum_{i\in G^z} \sum_{a \in \mathbb{A}} 1_{\{ X^i \in a\}}H^{z, i}_s \Lambda (ds \mid X^n_a,z) \quad (t \in \mathbb{R}_+)
\eeas
and 
\beas
\ddot{\mathcal{E}}^z_t
&=& 
\int_0^t 
\sum_{i\in G^z} \sum_{a \in \mathbb{A}} 1_{\{ X^i \in a\}}
H^{z, i}_s 
\big\{ \Lambda (ds \mid X^i,z) - \Lambda (ds \mid X^n_a,z) \big\}\quad (t \in \mathbb{R}_+),
\eeas
respectively, 
where  
$X^n_a$ is a random variable satisfying $X^n_a \in a$.
We have $\cale^z_t = \dot{\cale}^z_t + \ddot{\cale}^z_t$.

Define the event $\Omega^n$ by 
\begin{equation}\label{eq28}
\Omega^n
= 
\bigg\{\inf_{i\in G^1}\sum_{j\in G^0}1_{\{X^j\in a^i\}}Y^j_\tau>0\bigg\}.
\end{equation}
Condition $[A3]$ is equivalent to $P(\Omega^n) \to 1$ as $n \to \infty$.

\begin{lemma}\label{0410091235}
\bd
\im[(a)] 
On the event $\Omega^n$, 
\bea\label{0410081614}
\overline{\caly}^0_s 
&=& 
\overline{\caly}^1_s
\qquad( s\in[0,\tau]). 
\eea
In particular, $\overline{\caly}^0_0=\overline{\caly}^1_0=n_1$ on $\Omega^n$. 
\im[(b)] 
On the event $\Omega^n$, 
the functions $s\mapsto\overline{\caly}^z_s$ are non increasing for $z\in\{0,1\}$. 
\im[(c)] On the event $\Omega^n$, 
\beas
K^n_s 
&=&  
2^{-1/2}\big(\ol{\caly}^1_0\big)^{-1/2}\ol{\caly}^1_sW^n_s\qquad( s\in[0,\tau]). 
\eeas
\ed
\end{lemma}
\proof 
Since 
\beas
 1_{\{ X^{i'} \in a^i \}}\yeq 1_{\{ a^i=a^{i'} \}} \yeq 1_{\{ X^i \in a^{i'} \}}
\eeas
for $i, i' \in \bbG$ and
\beas 
\sum_{j\in G^0}1_{\{ X^j\in a^i \}}Y^j_{\tau} > 0\quad (\forall\omega \in \Omega^n, i\in G^1), 
\eeas
it follows that
\beas
\overline{\caly}^0_s 1_{\Omega^n}
&=& 
\sum_{i\in G^0} w^i_sY^i_s1_{\Omega^n} 
\nn\\&=&  
\sum_{i\in G^0}  
\left( \sum_{j\in G^0} 1_{\{ X^j \in a^i \}} Y^{j}_s\right)^{-}
\sum_{i'\in G^1}1_{\{ X^{i'} \in a^i \}} Y^{i'}_s
Y^i_s1_{\Omega^n} 
\nn\\&=& 
\sum_{i'\in G^1}\sum_{i\in G^0}  
\left( \sum_{j\in G^0} 1_{\{ X^j \in a^i \}} Y^{j}_s\right)^{-}
1_{\{ a^i=a^{i'} \}} Y^{i'}_sY^i_s1_{\Omega^n} 
\nn\\&=& 
\sum_{i'\in G^1}
\left( \sum_{j\in G^0} 1_{\{ X^j \in a^{i'} \}} Y^{j}_s\right)^{-}
\sum_{i\in G^0}  1_{\{ X^i \in a^{i'} \}} Y^i_sY^{i'}_s1_{\Omega^n} 
\nn\\&=& 
\sum_{i'\in G^1}Y^{i'}_s1_{\Omega^n}
\nn\\&=& 
\overline{\caly}^1_s 1_{\Omega^n}
\eeas
for all $s \in [0, \tau]$. 
Thus, we proved (a). 
We also have, on $\Omega^n$, 
$\overline{\caly}^0_s=\overline{\caly}^1_s=\sum_{i\in G^1}Y^i_s$, and the last function is nonincreasing, which shows (b). 
Since $\overline{\caly}^0_s = \overline{\caly}^1_s $ for $s\in[0,\tau]$ on $\Omega^n$ by (a), 
\beas
K^n_s 
&=&  
\big\{(\overline{\mathcal{Y}}^1_0 + \overline{\mathcal{Y}}^0_0)(\overline{\mathcal{Y}}^1_0\overline{\mathcal{Y}}^0_0)^-\big\}^{1/2}
(\overline{\mathcal{Y}}^1_s + \overline{\mathcal{Y}}^0_s)^-
\overline{\mathcal{Y}}^1_s \overline{\mathcal{Y}}^0_s W^n_s
\nn\\&=&
2^{-1/2}\big(\ol{\caly}^1_0\big)^{-1/2}\ol{\caly}^1_sW^n_s\qquad\text{on }\Omega^n, 
\eeas
that is (c). 
\qed\halflineskip
\begin{lemma}\label{0410081530}
Under $[H_0]$, 
$\dot{\cale}^1_t1_{\Omega^n}=\dot{\cale}^0_t1_{\Omega^n}$ for all $t\in[0,\tau]$ a.s. 
\end{lemma}
\proof
Fix any $s \in [0, \tau]$ and Borel set $B \in \bbB([0, \tau])$.
Since
\[
1_{\{ X^{i'} \in a^i \}}\yeq 1_{\{ a^i=a^{i'} \}} \yeq 1_{\{ X^i \in a^{i'} \}}
\]
for $i, i' \in \bbG$, we can calculate as follows on $\Omega^n$.
\begin{align}
&\sum_{i\in G^0} \sum_{a \in \mathbb{A}} 1_{\{ X^i \in a\}}w^i_sY^i_s \>\Lambda (B \mid X^n_a,0) \notag \\
=\; &\sum_{i\in G^0}\bigg(\sum_{j \in G^0}1_{\{X^j \in a^i\}}Y^j_s\bigg)^-
\sum_{i'\in G^1}1_{\{X^{i'}\in a^i\}}Y^{i'}_sY^i_s \>\Lambda (B \mid X^n_{a^i},0) \notag \\
=\; &\sum_{i\in G^0}\sum_{i'\in G^1}\bigg(\sum_{j\in G^0}1_{\{X^j \in a^i\}}Y^j_s\bigg)^- 1_{\{a^{i'} = a^i\}}Y^{i'}_s Y^i_s
\>\Lambda (B \mid X^n_{a^i},0) \notag \\
=\; &\sum_{i' \in G^1}\bigg(\sum_{j\in G^0}1_{\{X^j \in a^{i'}\}}Y^j_s\bigg)^- \sum_{i \in G^0}1_{\{X^i \in a^{i'}\}}Y^i_sY^{i'}_s
\Lambda (B \mid X^n_{a^{i'}},0) \notag \\
=\; &\sum_{i'\in G^1}Y^{i'}_s \Lambda (B \mid X^n_{a^{i'}},0) \notag \\
=\; &\sum_{i'\in G^1}\sum_{a \in \mathbb{A}}1_{\{X^{i'}\in a\}}w^{i'}_sY^{i'}_s
\Lambda (B \mid X^n_a,0).\label{0410081654} 
\end{align}
We use $[H_0]$, (\ref{0410081614}) of Lemma \ref{0410091235}, and (\ref{0410081654}), 
to show 
\beas &&
1_{\Omega^n} \sum_{i\in G^0} \sum_{a \in \mathbb{A}} 1_{\{ X^i \in a\}}H^{0, i}_s \Lambda (B \mid X^n_a,0) 
\nn\\&=& 
1_{\Omega^n} K^n_s\big(\ol{\caly}^0_s\big)^-\sum_{i\in G^0} \sum_{a \in \mathbb{A}} 1_{\{ X^i \in a\}}w^i_sY^i_s
\>\Lambda (B \mid X^n_a,1) 
\nn\\&=& 
1_{\Omega^n} K^n_s\big(\ol{\caly}^1_s\big)^-\sum_{i\in G^1} \sum_{a \in \mathbb{A}} 1_{\{ X^i \in a\}}w^i_sY^i_s
\>\Lambda (B \mid X^n_a,1). 
\eeas
Consequently, we proved the assertion of the lemma. 
\qed\halflineskip

\begin{lemma}\label{prop2}
If conditions $[A1]$, $[A3]$, $[A5]$, and $[A6]$ are satisfied, then
\beas
\sup_{t \in [0, \tau]} |\ddot{\mathcal{E}}^z_t| 
&=& 
O_p(\sqrt{n_1} d_n)\yeq o_p(1)
\eeas
as $n\to\infty$. 
\end{lemma}
\proof 
Let $z\in\{0,1\}$ and $t\in[0,\tau]$. 
From Lemma \ref{0410091235} (c), 
\beas
\big|H^{z, i}_s\big|
&=& 
\big|K^n_s (\overline{\caly}^z_s)^- w^i_sY^i_s\big|
\nn\\&=&
2^{-1/2}\big(\ol{\caly}^1_0\big)^{-1/2}|W^n_s|w^i_sY^i_s
\qquad\text{on }\Omega^n
\eeas
for $s \in [0, t]$ and $i\in\bbG$. 
%
%
Then, on $\Omega^n$, 
\beas
\big|\ddot{\mathcal{E}}^z_t\big|
&=& 
\bigg|\int_0^t 
\sum_{i\in G^z} \sum_{a \in \mathbb{A}} 1_{\{ X^i \in a\}}
H^{z, i}_s 
\big\{p(s\mid X^i,z) - p(s \mid X^n_a,z) \big\}\rho(ds)\bigg|
\nn\\&\leq& 
\int_0^t \sum_{i\in G^z} \sum_{a \in \mathbb{A}} 1_{\{ X^i \in a\}}|H^{z, i}_s|\rho(ds)Ld_n\qquad(\because [A1])
\nn\\&\leq& 
Ld_n\int_0^t \sum_{i\in G^z} \sum_{a \in \mathbb{A}} 1_{\{ X^i \in a\}}2^{-1/2}\big(\ol{\caly}^1_0\big)^{-1/2}|W^n_s|w^i_sY^i_s\rho(ds)
\nn\\&\leq& 
Ld_n\sup_{r\in[0,\tau]}|W^n_r|\>
\int_0^t \big(\ol{\caly}^1_0\big)^{-1/2}\sum_{i\in G^z} \sum_{a \in \mathbb{A}} 1_{\{ X^i \in a\}}w^i_sY^i_s\rho(ds)
\nn\\&\leq& 
Ld_n\sup_{r\in[0,\tau]}|W^n_r|
\int_0^\tau \big(\ol{\caly}^1_0\big)^{-1/2}\ol{\caly}^1_s\rho(ds)
\nn\\&\leq& 
n_1^{1/2}d_n\>L\sup_{r\in[0,\tau]}|W^n_r|\rho([0,\tau]), 
\eeas 
since $\big(\ol{\caly}^1_0\big)^{-1/2}\ol{\caly}_s\leq\big(\ol{\caly}^1_0\big)^{1/2}\leq n_1$. 
Remark that $s\mapsto\ol{\caly}^1_s$ is nonincreasing according to Lemma \ref{0410091235} (b). 
Consequently, we obtain the desired estimate since $P(\Omega^n) \to 1$ by $[A3]$ 
and $\sup_{r\in[0,\tau]}|W^n_r|=O_p(1)$ by $[A6]$. 
\qed\halflineskip

Next, we calculate $\langle \mathcal{M}^z \rangle$ to prove asymptotic normality of $\mathcal{M}^z$ 
by applying Theorem \ref{thm1}.

\begin{lemma}\label{0410111142} 
Suppose that $[A1]$, $[A2]$, $[A4]$ and $[A6]$ are fulfilled. Then
\bea\label{04101120145} 
\langle \mathcal{M}^0\rangle_t  &\to^p&0
\eea
and 
\bea\label{0410120146}
\langle \mathcal{M}^1\rangle_t 
&\to^p& 
\half\int_{[0,t]} (W^\infty_s)^2q_s\rho(ds)
\eea
as $n\to\infty$ for every $t\in[0,\tau]$. 
\end{lemma}
\proof 
Let $t\in[0,\tau]$. 
Since $[A3]$ implies $P(\Omega^n) \to 1$ as $n \to \infty$, it suffices to show that 
$\langle \mathcal{M}^0\rangle_t 1_{\Omega^n} \to 0$ and 
$\langle \mathcal{M}^1\rangle_t 1_{\Omega^n} \to \half\int_{[0,t]} (W^\infty_s)^2q_s\rho(ds)$. 
On $\Omega^n$, we have
\beas 
\langle \mathcal{M}^1\rangle_t 
&=& 
\sum_{i\in G^1}\int_0^t \big\{K^n_s(\overline{\mathcal{Y}}^1_s)^-w^i_s\big\}^2 d\langle M^i\rangle_s
\nn\\&=&
\sum_{i\in G^1}\int_0^t \big\{K^n_s(\overline{\mathcal{Y}}^1_s)^-w^i_s\big\}^2 Y^i_s \Lambda^i(ds \mid X^i, Z^i)
\nn\\&=&
\sum_{i\in G^1}\int_0^t 2^{-1}n_1^{-}(W^n_s)^2Y^i_s \Lambda^i(ds \mid X^i, Z^i)
\qquad(\because \text{Lemma }\ref{0410091235})
\nn\\&=&
2^{-1}\int_0^t (W^n_s)^2 n_1^{-}\sum_{i\in G^1}Y^i_s \>p(s\mid X^i, Z^i)\rho(ds).
\eeas
Then, the convergence (\ref{0410120146}) follows from $[A6]$ and $[A2]$.
Similarly, On $\Omega^n$, 
\beas 
\langle \mathcal{M}^0\rangle_t 
&=& 
\sum_{i\in G^0}\int_0^t \big\{K^n_s(\overline{\mathcal{Y}}^1_s)^-w^i_s\big\}^2 d\langle M^i\rangle_s
\nn\\&=&
\sum_{i\in G^0}\int_0^t \big\{K^n_s(\overline{\mathcal{Y}}^1_s)^-w^i_s\big\}^2 Y^i_s \Lambda^i(ds \mid X^i, Z^i)
\nn\\&=&
n_1^{-}\sum_{i\in G^0}\int_0^t 2^{-1}(W^n_s)^2Y^i_s (w^i_s)^2\Lambda^i(ds \mid X^i, Z^i)
\qquad(\because \text{Lemma }\ref{0410091235})
\nn\\&\leq&
\sup_{r\in[0,\tau]}(W^n_r)^2\cdot\sup_{r\in[0,\tau],j\in G^0}w^j_r\cdot
n_1^{-}\sum_{i\in G^0}\int_0^t 2^{-1}Y^i_s w^i_s\Lambda^i(ds \mid X^i, Z^i)
\nn\\&=&
\sup_{r\in[0,\tau]}(W^n_r)^2\cdot\sup_{r\in[0,\tau],j\in G^0}w^j_r\cdot
n_1^{-}\sum_{i\in G^0}\int_0^t 2^{-1}Y^i_s w^i_s\>p(s \mid X^i, Z^i)\rho(ds)
\nn\\&\leq&
\sup_{r\in[0,\tau]}(W^n_r)^2\cdot\sup_{r\in[0,\tau],j\in G^0}w^j_r\cdot
n_1^{-}\int_0^t 2^{-1}\sum_{i\in G^0}Y^i_s w^i_s\> \|p\|_\infty \rho(ds)
\nn\\&=&
\sup_{r\in[0,\tau]}(W^n_r)^2\cdot\sup_{r\in[0,\tau],j\in G^0}w^j_r\cdot
n_1^{-}\int_0^t 2^{-1}\ol{\caly}^1_s\> \|p\|_\infty \rho(ds)
\nn\\&\leq&
\sup_{r\in[0,\tau]}(W^n_r)^2\cdot\sup_{r\in[0,\tau],i\in G^0}w^i_r\cdot
2^{-1}\rho([0,\tau]) \|p\|_\infty. 
\eeas
By $[A1]$, $[A4]$ and $[A6]$, we obtain (\ref{0410120146}).
\qed\halflineskip

\begin{lemma}\label{lem4}
Suppose that $[C]$ and (\ref{0410130535}) of $[A1]$ are satisfied. Then
\bea\label{0512080945}
P\big(\Delta N^i_t = \Delta N^j_t \ge 1 \text{ for some }t \in [0, \tau]\big) &= &0
\eea
for $i\not=j$. 
\end{lemma}
\proof 
The process $[M^i,M^j]-\langle M^i,M^j\rangle\in\call$, the set of local martingales with the initial value $0$. 
Suppose that $i\not=j$. Then  
$R:=[M^i,M^j]\in\call$ by Proposition \ref{2209250620}. 
On the other hand, 
the local martingale $R$ is purely discontinuous because it has finite variation, i.e., 
the continuous martingale part of $R$ is null. 
Therefore, by Theorem 4.52 in \cite{jacod2003limit}, 
$R=\sum_{s\leq\cdot}\Delta M^i_s\Delta M^j_s=\sum_{s\leq\cdot}\Delta N^i_s\Delta N^j_s$, 
hence $0\leq R\leq1$, in particular, $R$ is a martingale. 
Thus, for any $t\in\bbR_+$, 
$0=E[R_0]=E[R_t]=E\big[\sum_{s\leq t}\Delta N^i_s\Delta N^j_s\big]$, which proves (\ref{0512080945}). 
\qed\halflineskip\halflineskip
%

%
\noindent
{\it Proof of Theorem \ref{0410120204}.}
We obtain the convergence 
$\sup_{t \in [0, \tau]} |\calm^0_t| \to^p 0$ as $n\to\infty$ 
by Lenglart's inequality and (\ref{04101120145}) of Lemma \ref{0410111142}.
Let $\epsilon > 0$.
Set $\cala_t = \sum_{s \le t} |\Delta \mathcal{M}^1_s |^2 1_{\{ |\Delta \mathcal{M}^1_s| > \epsilon \}}$ for $t \in [0, \tau]$.
We will show that $\cala_\tau^p \to^p 0$ as $n \to \infty$, where $\cala^p$ is the compensator of $\cala$.
Let $t \in [0, \tau]$.
Since $ A^i$ is continuous, we have
\begin{align*}
\cala_t 
&= \sum_{s \le t} \left( \sum_{i \in G^1} K^n_s (\overline{\mathcal{Y}}^1_s)^{-} \Delta M^i_s \right)^2
1_{\left\{\left| \sum_{j \in G^1} K^n_s (\overline{\mathcal{Y}}^1_s)^{-} \Delta M^j_s \right| > \epsilon\right\}} \\
&= \sum_{s \le t} \left( \sum_{i \in G^1} K^n_s (\overline{\mathcal{Y}}^1_s)^{-} \Delta N^i_s \right)^2
1_{\left\{\left| \sum_{j \in G^1} K^n_s (\overline{\mathcal{Y}}^1_s)^{-} \Delta N^j_s \right| > \epsilon\right\}}.
\end{align*}
By Lemma \ref{lem4},
\begin{align*}
&\sum_{s \le t} \left( \sum_{i \in G^1} K^n_s (\overline{\mathcal{Y}}^1_s)^{-} \Delta N^i_s \right)^2
1_{\left\{\left| \sum_{j \in G^1} K^n_s (\overline{\mathcal{Y}}^1_s)^{-} \Delta N^j_s \right| > \epsilon\right\}} \\
=\; &\sum_{s \le t} \sum_{i \in G^1} \left( K^n_s (\overline{\mathcal{Y}}^1_s)^{-} \Delta N^i_s \right)^2
1_{\left\{ \sum_{j \in G^1} \left| K^n_s (\overline{\mathcal{Y}}^1_s)^{-} \Delta N^j_s \right| > \epsilon\right\}} \\
=\; &\sum_{s \le t} \sum_{i \in G^1} \left( K^n_s (\overline{\mathcal{Y}}^1_s)^{-} \Delta N^i_s \right)^2
1_{\left\{ \left| K^n_s (\overline{\mathcal{Y}}^1_s)^{-} \Delta N^i_s \right| > \epsilon\right\}} \\
=\; &\sum_{s \le t} \sum_{i \in G^1} \left( K^n_s (\overline{\mathcal{Y}}^1_s)^{-} \right)^2
1_{\left\{ \left| K^n_s (\overline{\mathcal{Y}}^1_s)^{-} \right| > \epsilon\right\}} \Delta N^i_s \\
=\; &\sum_{i \in G^1} \int_0^t  \left( K^n_s (\overline{\mathcal{Y}}^1_s)^{-} \right)^2
1_{\left\{ \left| K^n_s (\overline{\mathcal{Y}}^1_s)^{-} \right| > \epsilon\right\}} dN^i_s.
\end{align*}
Therefore, it holds that
\[
\cala^p_t = \sum_{i \in G^1} \int_0^t  \left( K^n_s (\overline{\mathcal{Y}}^1_s)^{-} \right)^2
1_{\left\{ \left| K^n_s (\overline{\mathcal{Y}}^1_s)^{-} \right| > \epsilon\right\}} dA^i_s.
\]
By Lemma \ref{0410091235} (c), we have
\begin{align*}
\cala^p_\tau 1_{\Omega^n}
&\le \sum_{i \in G^1} \int_0^\tau (2n_1)^-(W^n_s)^2 1_{\{(2n_1)^- (W^n_s)^2 > \epsilon^2 \}} dA^i_s 1_{\Omega^n} \\
&\le \sum_{i \in G^1} \int_0^\tau (4n_1^2 \epsilon^2)^- (W^n_s)^4 1_{\{(2n_1)^- (W^n_s)^2 > \epsilon^2 \}} dA^i_s 1_{\Omega^n} \\
&\le \sup_{s \in [0, \tau]} (W^n_s)^4 \cdot (4n_1^2 \epsilon^2)^- \sum_{i \in G^1} A^i_\tau.
\end{align*}
Conditions $[A5]$ and $[A6]$ imply that $\sup_{s \in [0, \tau]} (W^n_s)^4 \cdot n_1^- \to^p 0$.
Moreover, $\left\{n_1^- \sum_{i \in G^1} A^i_\tau\right\}_{n \in \bbN}$ is tight since 
\[
E\left[ n_1^- \sum_{i \in G^1} A^i_\tau \right] 
= E\left[ n_1^- \sum_{i \in G^1} E[A^i_\tau \mid \calf_0] \right]
= E\left[ n_1^- \sum_{i \in G^1} E[N^i_\tau \mid \calf_0] \right]
\le 1.
\]
Thus, we have $\cala^p_\tau 1_{\Omega^n} \to^p 0$.
By $[A3]$, we obtain $\cala^p_\tau \to^p 0$.
From this convergence together with (\ref{0410120146}) of Lemma \ref{0410111142}, we apply 
Rebolledo's martingale central limit theorem (Theorem \ref{thm1}) 
to obtain the convergence $\calm^1\to^d\bbW$. 
Since the sequences $\calm^0$ and $\cale^1 - \cale^0$ are C-tight by Lemmas \ref{0410081530} and \ref{prop2}, 
the decomposition (\ref{0410121858}) derives the convergence $\calw\to^d\bbW$ in $D([0,\tau])$. 
\qed

\section{Proof of Proposition \ref{prop4}}\label{section2}
Since $W^\infty_s$ is bounded, $[A2]$ implies
\bea
&&\left| (2n_1)^- \sum_{i \in G^1} \int_0^\tau (W^\infty_s)^2\,dA^i_s - \text{Var}[\bbW_\tau] \right| \notag 
\nn\\
&\le&
\frac{1}{2} \int_0^\tau (W^\infty_s)^2 
\left| n_1^- \sum_{i \in G^1} Y^i_s p(s \mid X^i, Z^i) - q_s \right|\,\rho(ds) \notag 
\nn\\
&\le&
\frac{1}{2} \sup_{s \in [0, \tau]} (W^\infty_s)^2 \cdot
\sup_{s \in [0, \tau]} \left| n_1^- \sum_{i \in G^1} Y^i_s p(s \mid X^i, Z^i) - q_s \right| \cdot \rho([0, \tau]) \notag 
\nn\\
&\to^p&
0. \label{eq21}
\eea
We have
\begin{align}
&\left| (2n_1)^- \sum_{i \in G^1} \int_0^\tau (W^n_s)^2\,dA^i_s 
- (2n_1)^- \sum_{i \in G^1} \int_0^\tau (W^\infty_s)^2\,dA^i_s \right| \notag \\
\le\; &\sup_{s \in [0, \tau]} |(W^n_s)^2 - (W^\infty_s)^2| \cdot (2n_1)^- \sum_{i \in G^1} A^i_\tau \notag \\
\to^p\; &0, \label{eq22}
\end{align}
where we used $[A6]$ and 
the estimate $A^i_\tau\leq\|p\|_\infty\rho([0,\tau])$.
Furthermore, using Lenglart's inequality to the locally square-integrable martingale, 
we see 
\bea\label{0512061557}
\calv_\tau - (2n_1)^- \sum_{i \in G^1} \int_0^\tau (W^n_s)^2\,dA^i_s
&=& 
(2n_1)^- \sum_{i \in G^1} \int_0^\tau (W^n_s)^2\,dM^i_s 
\>\to^p\>
0
\eea
since
\beas 
\left\langle n_1^- \sum_{i \in G^1} \int_0^\cdot (W^n_s)^2\,dM^i_s \right\rangle_\tau
&=&
(n_1^-)^2 \sum_{i \in G^1} \int_0^\tau (W^n_s)^4 \,dA^i_s\\
&\yleq&
n_1^-\sup_{s\in[0,\tau]}|W^n_s|^4 \|p\|_\infty \rho([0, \tau])\\
&\>\to^p\>&
0. 
\eeas
%
It follows from equations (\ref{eq21}), (\ref{eq22}) and (\ref{0512061557}) that $\calv_\tau \to^p \text{Var}[\bbW_\tau]$.
\qed

\section{Proof of Theorem \ref{0510271447}}\label{0510271450}
We have shown that
\bea\label{eq24}
\dot{\cale}^1_\tau
&=&
 \int_0^\tau \sum_{i \in G^1} \sum_{a \in \mathbb{A}}1_{\{X^i \in a\}} H^{1, i}_s \Lambda(ds \mid X^n_a, 1) \notag 
 \nn\\&=& 
\int_0^\tau \sum_{i \in G^1} \sum_{a \in \mathbb{A}}1_{\{X^i \in a\}} H^{1, i}_s \Lambda(ds \mid X^i, 1) + o_p(1)
\eea
in Lemma \ref{prop2}. 
By (\ref{0410081654}) of Lemma \ref{0410081530} and the proof of Lemma \ref{prop2}, it holds that
\bea\label{eq25}
\dot{\cale}^0_t
&=& 
\int_0^\tau \sum_{i \in G^0} \sum_{a \in \mathbb{A}} 1_{\{X^i \in a\}} H^{0, i}_s \Lambda(ds \mid X^n_a, 0) \notag 
\nn\\&=& 
\int_0^\tau \sum_{i \in G^1} \sum_{a \in \mathbb{A}}1_{\{X^i \in a\}} H^{1, i}_s \Lambda(ds \mid X^n_a, 0) + o_p(1) \notag
\nn\\&=& 
\int_0^\tau \sum_{i \in G^1} \sum_{a \in \mathbb{A}}1_{\{X^i \in a\}} H^{1, i}_s \Lambda(ds \mid X^i, 0) + o_p(1).
\eea
Equations (\ref{eq24}), (\ref{eq25}) and Lemma \ref{0410091235} implies
\beas
&&\dot{\mathcal{E}}^1_\tau - \dot{\mathcal{E}}^0_\tau\\
&=&  
\int_0^\tau \sum_{i \in G^1} \sum_{a \in \mathbb{A}} 1_{\{X^i \in a\}} H^{1, i}_s
\{ \Lambda(ds \mid X^i, 1) - \Lambda(ds \mid X^i, 0) \} + o_p(1) 
\\&=& 
\int_0^\tau \sum_{i \in G^1} K^n_s (\overline{\mathcal{Y}}^1_s)^- Y^i_s
\{ \Lambda(ds \mid X^i, 1) - \Lambda(ds \mid X^i, 0) \} + o_p(1)
\\&=& 
\int_0^\tau W^n_s \, 2^{-1/2} n_1^{-1/2} \sum_{i \in G^1} 
Y^i_s\{ p(s \mid X^i, 1) - p(s \mid X^i, 0) \}\,\rho(ds) + o_p(1).
\eeas
Therefore, by $[A7]$ and $[H1]$, we have
\beas
&&(n_1^-)^{1/2}(\dot{\mathcal{E}}^1_\tau - \dot{\mathcal{E}}^0_\tau)\\
&=& 
\int_0^\tau W^n_s\,2^{-1/2} n_1^- \sum_{i \in G^1} 
Y^i_s \{ p(s \mid X^i, 1) - p(s \mid X^i, 0) \} \,\rho(ds) \,1_{\Omega^n} + o_p(1) 
\\&\le& 
\sup_{s \in B, x \in \mathcal{X}} W^n_s \{p(s \mid x, 1) - p(s \mid x, 0)\} 
\cdot n_1^- \sum_{i \in G^1} Y^i_\tau \cdot \rho(B) + o_p(1) 
\\&\to^p& 
\sup_{s \in B, x \in \mathcal{X}} W^\infty_s \{p(s \mid x, 1) - p(s \mid x, 0)\} 
\cdot \widetilde{q}_\tau \cdot \rho(B) 
\\&<& 
0.
\eeas
Since $\ddot{\mathcal{E}}^z_\tau = o_p(1)$ ($z = 0, 1$) by Lemma \ref{prop2} and
$\mathcal{M}^z_\tau = O_p(1)$ ($z = 0, 1$) by the proof of Theorem \ref{0410120204},
it holds that $P(\mathcal{W}_\tau < -R) \to 1$ as $n \to \infty$ for all $R > 0$.
\qed

\section{Proof of Theorem \ref{0510271457}}\label{0510271458}
By Theorem \ref{0410120204}, Proposition \ref{prop4} and Theorem \ref{0510271447}, 
if $[C]$ and $[A1]$-$[A7]$ are fulfilled, then $(1)$-$(3)$ hold.
Thus, it is enough to show $[A2]$-$[A5]$ and $[A7]$.

First, we show $[A3]$ and $[A4]$ from $[B1]$-$[B3]$ and $[B5]$.
Let $\delta, \delta' > 0$ be constants satisfying $\delta' < \delta < (1 - \overline{d}\theta) \wedge (1 - \beta)$.
For $n = 1, 2, \ldots$, we define the events $\Omega^{n, 1}$ and $\Omega^{n, 0}$ as
\begin{align*}
\Omega^{n, 1} 
&= \left\{ \max_{a \in \mathbb{A}} \sum_{i \in \mathbb{G}^1} 1_{\{X^i \in a\}} < n^{1-\overline{d}\theta - \delta} \right\},\\
\Omega^{n, 0} 
&= \left\{ \min_{a \in \mathbb{A}} \sum_{i \in \mathbb{G}^0} 1_{\{X^i \in a\}}Y^i_\tau > n^{1- \overline{d}\theta - \delta'} 
\right\}, 
\end{align*}
respectively.
Condition $[A3]$ holds if $P(\Omega^{n, 0}) \to 1$ ($n \to \infty$).
On the event $\Omega^{n, 0} \cap \Omega^{n, 1}$, we have
\[
\sup_{t \in [0, \tau], i \in G^0} w^i_t
\le \max_{a \in \mathbb{A}}\left\{ \left( \sum_{j \in \mathbb{G}^0} 1_{\{X^j \in a\}} Y^j_\tau \right)^{-1} 
\sum_{j' \in \mathbb{G}^1}1_{\{X^{j'} \in a\}}\right\}
 < n^{-(\delta - \delta')}.
\]
Therefore, $[A4]$ holds if $P(\Omega^{n, 0}) \to 1$ and $P(\Omega^{n, 1}) \to 1$ ($n \to \infty$).

First, we show $P(\Omega^{n, 1}) \to 1$.
Let $p^{n, 1}(a) = P(X^1 \in a, Z^1 = 1)$ for each $n \in \mathbb{N}$ and $a \in \mathbb{A}$.
It holds that
\begin{align*}
p^{n, 1}(a) 
&= P(X^1 \in a \mid Z^1 = 1) P(Z^1 = 1) \\
&= \int_a f(x \mid 1) \,\nu(dx) \cdot P(Z^1 = 1) \\
&\le \esssup_{x \in \mathcal{X}} f(x \mid 1) \cdot \nu(a)P(Z^1 = 1).
\end{align*}
Thus, by $[B2]$, $[B3]$ and $[B5]$, we have $\max_{a \in \mathbb{A}} p^{n, 1}(a) \lesssim n^{\beta - 1 - \overline{d}\theta}$.
Let $(\xi^{i})_{i = 1}^n$ be an i.i.d.~random variables that follow the standard uniform distribution $U(0, 1)$.
Set $k^{n, 1} = \lceil n^{1 - \overline{d}\theta - \delta} \rceil$.
Since $\sum_{i \in \mathbb{G}^1} 1_{\{X^i \in a\}}$ follows the binomial distribution $B(n, p^{n, 1}(a))$, 
we have
\begin{align*}
P\left( \sum_{i \in \mathbb{G}^1} 1_{\{X^i \in a\}} \ge n^{1 - \overline{d}\theta - \delta} \right) 
&= P\left( \sum_{i = 1}^n 1_{\{\xi^i \le p^{n, 1}(a)\}} \ge \lceil n^{1 - \overline{d}\theta - \delta} \rceil \right) \\
&= P(\xi^{(k^{n, 1})} \le p^{n, 1}(a)), 
\end{align*}
where $\xi^{(k)}$ is the $k$-th order statistic of $(\xi^i)_{i = 1}^n$.
The order statistic $\xi^{(k^{n, 1})}$ follows the beta distribution $\text{Beta}(k^{n, 1}, n - k^{n, 1} + 1)$.
Thus, by Stirling's formula, 
\begin{align*}
P(\xi^{(k^{n, 1})} \le p^{n, 1}(a))
&= \frac{n!}{(k^{n, 1} - 1)!(n - k^{n, 1})!} \int_0^{p^{n, 1}(a)} x^{k^{n, 1} - 1} (1-x)^{n - k^{n, 1}}\,dx \\
&\le \frac{n^{k^{n, 1}}}{(k^{n, 1} - 1)!} \int_0^{p^{n, 1}(a)} x^{k^{n, 1} - 1}\,dx \\
&= \frac{(np^{n, 1}(a))^{k^{n, 1}}}{k^{n, 1}!} \\
&\le \frac{ (np^{n, 1}(a))^{k^{n, 1}} e^{k^{n, 1}} }{ {k^{n, 1}}^{k^{n, 1}} (2\pi k^{n, 1})^{1/2} }\\
&\le \left(\frac{ enp^{n, 1}(a) }{ k^{n, 1} }\right)^{k^{n, 1}}.
\end{align*}
Since 
$\max_{a \in \mathbb{A}} np^{n, 1}(a) \lesssim n^{\beta - \overline{d}\theta} 
= n^{1 - \overline{d}\theta - (1 - \beta)}
= o(n^{1 - \overline{d}\theta - \delta})
= o(k^{n, 1})$, 
it holds that 
\[
\max_{a \in \mathbb{A}} \left(\frac{ enp^{n, 1}(a) }{ k^{n, 1} }\right)^{k^{n, 1}} 
\le e^{-k^{n, 1}} 
\le e^{-n^{1 - \overline{d}\theta - \delta}}
\]
for large enough $n \in \mathbb{N}$.
Hence, we obtain
\[
\max_{a \in \mathbb{A}} P\left( \sum_{i \in \mathbb{G}^1} 1_{\{X^i \in a\}} \ge n^{1 - \overline{d}\theta - \delta} \right)
\lesssim e^{-n^{1 - \overline{d}\theta - \delta}}.
\]
Therefore, by $[B3]$ and $[B5]$, 
\begin{align*}
P((\Omega^{n, 1})^c)
&= P\left( \max_{a \in \mathbb{A}} \sum_{i \in \mathbb{G}^1} 1_{\{X^i \in a\}} \ge n^{1-\overline{d}\theta - \delta} \right) \\
&\le \sum_{a \in \mathbb{A}} P\left( \sum_{i \in \mathbb{G}^1} 1_{\{X^i \in a\}} \ge n^{1-\overline{d}\theta - \delta} \right) \\
&\lesssim \#\mathbb{A} \cdot e^{-n^{1 - \overline{d}\theta - \delta}} \\
&\le \left( \min_{a \in \mathbb{A}} \nu(a) \right)^{-1} \nu(\mathcal{X}) \cdot e^{-n^{1 - \overline{d}\theta - \delta}} \\
&\lesssim n^{\overline{d}\theta} e^{-n^{1 - \overline{d}\theta - \delta}} \\
&\to 0.
\end{align*}

Next, we show $P(\Omega^{n, 0}) \to 1$.
Let $p^{n, 0}(a) = P(X^1 \in a, Z^1 = 0, \widetilde{T}^1 \ge \tau)$ for each $n \in \mathbb{N}$ and $a \in \mathbb{A}$.
From $[B2]$, $[B3]$, $[B5]$ and 
\begin{align*}
p^{n, 0}(a) 
&= \int_{a \times \{0\}} P(\widetilde{T}^1 \ge \tau \mid X^1 = x, Z^1 = z) \,P^{(X^1, Z^1)}(dxdz) \\
&= \int_a P(\widetilde{T}^1 \ge \tau \mid X^1 = x, Z^1 = 0) f(x \mid 0) \,\nu(dx) \cdot P(Z^1 = 0) \\
&\ge \essinf_{x \in \mathcal{X}} \left( P(\widetilde{T}^1 \ge \tau \mid X^1 = x, Z^1 = 0) f(x \mid 0) \right)
\nu(a) P(Z^1 = 0), 
\end{align*}
it follows that $\min_{a \in \mathbb{A}} p^{n, 0}(a) \gtrsim n^{-\overline{d}\theta}$.
Set $k^{n, 0} = \lfloor n^{1 - \overline{d}\theta - \delta'} \rfloor + 1$.
By the same way as in the proof of $P(\Omega^{n, 1}) \to 1$, we can calculate as follows.
\begin{align*}
&P\left( \sum_{i \in \mathbb{G}^0} 1_{\{X^i \in a\}}Y^i_\tau \le n^{1- \overline{d}\theta - \delta'} \right) \\
=\; &P\left( \sum_{i = 1}^n 1_{\{\xi^i \le p^{n, 0}(a)\}} \le \lfloor n^{1 - \overline{d}\theta - \delta'} \rfloor \right) \\
=\; &P(\xi^{(k^{n, 0})} > p^{n, 0}(a)) \\
=\; &\frac{n!}{(k^{n, 0} - 1)! (n - k^{n, 0})!} \int_{p^{n, 0}(a)}^1 x^{k^{n, 0} - 1} (1 - x)^{n - k^{n, 0}} \,dx \\
\le\; &\frac{n!}{(k^{n, 0} - 1)! (n - k^{n, 0})!} \int_{p^{n, 0}(a)}^1 (1 - x)^{n - k^{n, 0}} \,dx \\
=\; &\frac{n!}{(k^{n, 0} - 1)! (n - k^{n, 0} + 1)!} (1 - p^{n, 0}(a))^{n - k^{n, 0} + 1}.
\end{align*}
By the stirling's formula, for large enough $n \in \mathbb{N}$, we have
\begin{align}
&\log P\left( \sum_{i \in \mathbb{G}^0} 1_{\{X^i \in a\}}Y^i_\tau \le n^{1- \overline{d}\theta - \delta'} \right) \notag \\
\le\; &\log n! - \log (k^{n, 0} - 1)! - \log (n - k^{n, 0} + 1)! + (n - k^{n, 0} + 1) \log (1 - p^{n, 0}(a)) \notag \\
\le\; &(n + 1/2) \log n - n 
- (k^{n, 0} - 1/2) \log (k^{n, 0} - 1) + k^{n, 0} - 1 \notag \\
&\quad - (n - k^{n, 0} + 3/2) \log (n - k^{n, 0} + 1) + n - k^{n, 0} + 1 \notag \\
&\quad + (n - k^{n, 0} + 1) \log (1 - p^{n, 0}(a)) + C \notag \\
=\; &(n + 1/2) \log n - (k^{n, 0} - 1/2) \log (k^{n, 0} - 1) \notag \\
&\quad - (n - k^{n, 0} + 3/2) \left(\log n + \log \left(1 - \frac{k^{n, 0} - 1}{n}\right) \right) \notag \\
&\quad + (n - k^{n, 0} + 1) \log (1 - p^{n, 0}(a)) + C \notag \\
=\; &(n + 1/2) \log n - (k^{n, 0} - 1/2) \log (k^{n, 0} - 1) \notag \\
&\quad - (n - k^{n, 0} + 3/2) \left(\log n - \frac{k^{n, 0} - 1}{n} + o\left( \frac{k^{n, 0} - 1}{n} \right)\right) \notag \\
&\quad - (n - k^{n, 0} + 1) p^{n, 0}(a) + C \notag \\
=\; &1/2 \cdot \log n - (k^{n, 0} - 1/2) \log (k^{n, 0} - 1) + (k^{n, 0} - 3/2) \log n \notag \\
&\quad + (n - k^{n, 0} + 3/2) \cdot \left( \frac{k^{n, 0} + 1}{n} + o\left( \frac{k^{n, 0} - 1}{n} \right) \right)
- (n - k^{n, 0} + 1) p^{n, 0}(a) + C \notag \\
\le\; &2k^{n, 0} \log n - np^{n, 0}(a)/2 + C, \label{eq49}
\end{align}
where $C$ is some constant.
Since $\min_{a \in \mathbb{A}} np^{n, 0}(a) \gtrsim n^{1 - \overline{d}\theta}$ 
and $k^{n, 0} \asymp n^{1 - \overline{d}\theta - \delta'}$, 
the principal term in (\ref{eq49}) is $-np^{n, 0}(a)/2$.
Therefore, there exists $b > 0$ such that 
\[
\max_{a \in \mathbb{A}} P\left( \sum_{i \in \mathbb{G}^0} 1_{\{X^i \in a\}}Y^i_\tau \le n^{1- \overline{d}\theta - \delta'} \right)
\le \exp\left( - \frac{1}{4}\min_{a \in \mathbb{A}} np^{n, 0}(a) \right)
\le e^{-b n^{1 - \overline{d}\theta}}
\]
for large enough $n \in \mathbb{N}$.
Thus, we have
\begin{align*}
P((\Omega^{n, 0})^c)
&= P\left( \min_{a \in \mathbb{A}} \sum_{i \in \mathbb{G}^0} 1_{\{X^i \in a\}} Y^i_\tau \le n^{1 - \overline{d}\theta - \delta'} \right) \\
&\le \sum_{a \in \mathbb{A}} P\left( \sum_{i \in \mathbb{G}^0} 1_{\{X^i \in a\}} Y^i_\tau 
\le n^{1 - \overline{d}\theta - \delta'} \right) \\
&\le \# \mathbb{A} \cdot e^{-bn^{1 - \overline{d}\theta}} \\
&\lesssim n^{\overline{d}\theta} e^{-bn^{1 - \overline{d}\theta}} \\
&\to 0.
\end{align*}
The proofs of $[A3]$ and $[A4]$ are now complete.

We recall that $V^i = (X^i, Z^i, T^i, U^i)$ ($i \in \bbG$).
Let $g \colon \bbR^d \times \{0, 1\} \times (0, \infty) \times (0, \infty) \to \bbR$ be a measurable function 
satisfying $E[|g(V^1)|^2 \mid Z^1 = 1] < \infty$.
Note that $E[|g(V^1)|^2 \mid Z^1 = 1]$ does not depend on $n$ 
since the conditional distribution $\call \{V^1 \mid Z^1 = 1\}$ does not depend on $n$.
Suppose that $[B1]$-$[B5]$ hold.
In what follows, we consider large enough $n \in \bbN$ such that $P(Z^1 = 1) > 0$.
We prove the weak law of large numbers
\begin{equation}\label{eq33}
\frac{1}{nP(Z^1 = 1)}\sum_{i \in G^1} g(V^i) \to^p E[g(V^1) 1_{\{X^1 \in \calx\}} \mid Z^1 = 1]
\end{equation}
and
\begin{equation}\label{eq36}
n_1^-\sum_{i \in G^1} g(V^i) \to^p E[g(V^1) \mid X^1 \in \calx, Z^1 = 1].
\end{equation}
On the event $\Omega^{n, 0}$, for all $a \in \bbA$, there exists at least one individual $i \in \bbG^0$ such that $X^i \in a$.
Therefore, we have 
\begin{align}
G^1 
&= \{ i \in \bbG^1 \mid \text{$X^i \in a$ and $X^j \in a$ for some $j \in \bbG^0$ and $a \in \bbA$} \} \notag \\
&= \{ i \in \bbG^1 \mid X^i \in \calx \} \quad \text{on $\Omega^{n, 0}$}. \label{eq34}
\end{align}
We have shown that if $[B1]$-$[B3]$ and $[B5]$ holds, then $P(\Omega^{n, 0}) \to 1$.
From $P(\Omega^{n, 0}) \to 1$ and (\ref{eq34}), it follows that
\begin{align*}
\frac{1}{nP(Z^1 = 1)}\sum_{i \in G^1} g(V^i)
&= \frac{1}{nP(Z^1 = 1)}\sum_{i \in G^1} g(V^i) 1_{\Omega^{n, 0}} + o_p(1) \\
&= \frac{1}{nP(Z^1 = 1)}\sum_{i \in \bbG^1} g(V^i) 1_{\{X^i \in \calx\}}1_{\Omega^{n, 0}} + o_p(1) \\
&= \frac{1}{nP(Z^1 = 1)}\sum_{i \in \bbG^1} g(V^i) 1_{\{X^i \in \calx\}} + o_p(1).
\end{align*}
Thus, we need to show that 
\begin{equation}\label{eq35}
\frac{1}{nP(Z^1 = 1)}\sum_{i \in \bbG^1} g(V^i) 1_{\{X^i \in \calx\}} \to^p E[g(V^1) 1_{\{X^1 \in \calx\}} \mid Z^1 = 1]
\end{equation}
for (\ref{eq33}).
Note that the strong law of large numbers may not hold in (\ref{eq35}) 
since the distribution of $(X^i, Z^i, T^i, U^i)$ depends on $n$.
We have
\begin{align*}
E\left[ \frac{1}{nP(Z^1 = 1)} \sum_{i \in \bbG^1} g(V^i)1_{\{X^i \in \calx\}} \right] 
&= \frac{1}{nP(Z^1 = 1)}  E \left[ \sum_{i \in \bbG} g(V^i) 1_{\{X^i \in \calx, Z^i = 1\}} \right] \\
&= \frac{1}{P(Z^1 = 1)} E[ g(V^1) 1_{\{X^1 \in \calx, Z^1 = 1\}} ] \\
&= E[g(V^1) 1_{\{X^1 \in \calx\}} \mid Z^1 = 1]
\end{align*}
and 
\begin{align*}
\text{Var}\left[ \frac{1}{nP(Z^1 = 1)} \sum_{i \in \bbG^1} g(V^i) 1_{\{X^i \in \calx\}} \right] 
&= \frac{1}{n^2P(Z^1 = 1)^2} \text{Var}\left[ \sum_{i \in \bbG} g(V^i) 1_{\{X^i \in \calx, Z^i = 1\}}\right] \\
&= \frac{1}{nP(Z^1 = 1)^2} \text{Var}[ g(V^1) 1_{\{X^1 \in \calx, Z^1 = 1\}} ] \\
&\le \frac{1}{nP(Z^1 = 1)^2} E[(g(V^1))^2 1_{\{X^1 \in \calx, Z^1 = 1\}}] \\
&= \frac{1}{nP(Z^1 = 1)} E[(g(V^1))^2 1_{\{X^1 \in \calx\}} \mid Z^1 = 1] \\
&= O_p(n^{-\beta})
\end{align*}
by $[B2]$.
Therefore, we have 
\[
\frac{1}{nP(Z^1 = 1)}\sum_{i \in \bbG^1} g(V^i) 1_{\{X^i \in \calx\}} \to E[g(V^1) 1_{\{X^1 \in \calx\}} \mid Z^1 = 1] 
\quad \text{in $L^2$}.
\]
Thus, we have (\ref{eq35}) and we obtain (\ref{eq33}).
By applying (\ref{eq33}) to $g = 1$, we have
\begin{equation}\label{eq37}
\frac{n_1}{nP(Z^1 = 1)} \to^p P(X^1 \in \calx \mid Z^1 = 1) > 0,
\end{equation}
where $P(X^1 \in \calx \mid Z^1 = 1) > 0$ is derived from $[B4]$.
It follows from (\ref{eq33}) and (\ref{eq37}) that
\[
n_1^-\sum_{i \in G^1} g(V^i) 
\to^p \frac{E[g(V^1) 1_{\{X^1 \in \calx\}} \mid Z^1 = 1]}{P(X^1 \in \calx \mid Z^1 = 1)}
\]
and (\ref{eq36}) is proved.

We derive $[A5]$ from $[B1]$-$[B5]$.
By (\ref{eq37}) and $[B2]$, it holds that $n_1 = O_p(n^{\beta})$ and $n_1^{-} = O_p(n^{-\beta})$.
Therefore, by $[B5]$ and $[B1]$, we have $n_1^{-} = o_p(1)$ and $n_1d_n^2 = O_p(n^{\beta-2\theta}) = o_p(1)$, 
which implies $[A5]$.

Condition $[A7]$ follows from $[B1]$-$[B5]$.
Indeed, when $g(x, z, t, u) = 1_{\{t \wedge u \ge \tau\}}$ ($x \in \bbR^d, z \in \{0, 1\}, t, u \in (0, \infty)$),
the weak law of large numbers (\ref{eq36}) gives
\[
n_1^- \sum_{i \in G^1} Y^i_\tau \to^p \widetilde{q}_\tau \coloneqq E[Y^1_\tau \mid X^1 \in \calx, Z^1 = 1].
\]
By [B4], we have 
\[
\widetilde{q}_\tau = \frac{P(\widetilde{T}^1 \ge \tau, X^1 \in \calx \mid Z^1 = 1)}{P(X \in \calx \mid Z^1 = 1)} > 0.
\]



Finally, we show $[A2]$ under $[A1']$ and $[B1]$-$[B5]$.
The sequence of \cadlag processes 
$\big(\bbV^n_s \coloneqq n_1^-\sum_{i\in G^1}Y^i_{s+}p(s \mid X^i,Z^i)\big)_{n\in\bbN}$ is tight in $D([0, \tau])$.
The proof is as follows.

For every $\alpha \in D([0, \tau])$ and $\delta > 0$, we define
\[
w(\alpha, \delta) = \sup_{s, t \in [0, \tau],\, |s-t| \le \delta} |\alpha_s - \alpha_t|.
\]
Condition $[A1']$ implies
\[
\sup_{n \in \bbN, s \in [0, \tau]} |\bbV^n_s| 
\le \|p\|_\infty = \sup_{s \in [0, \tau], x \in \calx, z \in \{0, 1\}} p(s \mid x, z)
< \infty. 
\]
Hence, it is enough to show that 
\begin{equation}\label{eq46}
\lim_{\delta \downarrow 0} \limsup_{n \to \infty} P(w(\bbV^n, \delta) \ge \epsilon) = 0 \quad \text{for all $\epsilon > 0$}
\end{equation}
to obtain the (C-)tightness of $(\bbV^n)$.
Let $\delta > 0$ and $n \in \bbN$.
If $s, t \in [0, \tau]$ satisfy $s \le t$ and $t - s \le \delta$, then
\begin{align*}
|\bbV^n_s - \bbV^n_t|
&\le \left| n_1^- \sum_{i \in G^1} Y^i_{s+}p(s \mid X^i, Z^i) - n_1^- \sum_{i \in G^1} Y^i_{t+}p(s \mid X^i, Z^i) \right| \\
&\quad +  \left| n_1^- \sum_{i \in G^1} Y^i_{t+}p(s \mid X^i, Z^i) - n_1^- \sum_{i \in G^1} Y^i_{t+}p(t \mid X^i, Z^i) \right| \\
&\le n_1^- \sum_{i \in G^1} (Y^i_{s+} - Y^i_{t+}) \|p\|_\infty
+ n_1^- \sum_{i \in G^1} |(p(s \mid X^i, Z^i) - p(t \mid X^i, Z^i))| \\
&\le \|p\|_\infty \left( n_1^- \sum_{i \in G^1} Y^i_{s+} - n_1^- \sum_{i \in G^1} Y^i_{t+} \right) 
+ 
\sup_{x \in \calx} |p(s \mid x, 1) - p(t \mid x, 1)| \\
&\le \|p\|_\infty w\left( n_1^-\sum_{i \in G^1}Y^i_{+}, \delta \right) 
+ \sup_{s, t \in [0, \tau],\, |s-t| \le \delta,\, x \in \calx} |p(s \mid x, 1) - p(t \mid x, 1)|, 
\end{align*}
where $Y^i_+$ denotes the process $(Y^i_{s+})_{s \in [0, \tau]}$.
By taking the supremum over $s$ and $t$, we have
\[
w(\bbV^n, \delta) 
\le \|p\|_\infty w\left( n_1^-\sum_{i \in G^1}Y^i_{+}, \delta \right) 
+ \sup_{s, t \in [0, \tau],\, |s-t| \le \delta,\, x \in \calx} |p(s \mid x, 1) - p(t \mid x, 1)|.
\]
Therefore, by Proposition 3.26 in Chapter 6 of \cite{jacod2003limit}, we need to show that
\begin{equation}\label{eq47}
\lim_{\delta \downarrow 0} \limsup_{n \to \infty} 
P\left( w\left( n_1^-\sum_{i \in G^1}Y^i_{+}, \delta \right) \ge \epsilon \right) = 0 \quad \text{for all $\epsilon > 0$}
\end{equation}
and
\begin{equation}\label{eq48}
\lim_{\delta \downarrow 0} \sup_{s, t \in [0, \tau],\, |s-t| \le \delta,\, x \in \calx} |p(s \mid x, 1) - p(t \mid x, 1)| = 0
\end{equation}
to obtain (\ref{eq46}).

By $[A1']$, there exists $L > 0$ such that for all $s, t \in [0, \tau]$ and $x, y \in \calx$
\begin{align*}
|p(s \mid x, 1) - p(t \mid y, 1)| 
&\le |p(s \mid x, 1) - p(t \mid x, 1)| + |p(t \mid x, 1) - p(t \mid y, 1)| \\
&\le |p(s \mid x, 1) - p(t \mid x, 1)| + L|x - y|.
\end{align*}
From this and $[A1']$, the function $(s, x) \mapsto p(s \mid x, 1)$ is continuous on $[0, \tau] \times \calx$.
Since $\calx$ is compact, $p(s \mid x, 1)$ is uniformly continuous on $[0, \tau] \times \calx$.
Therefore, (\ref{eq48}) holds.

For each $s \in [0, \tau]$, applying (\ref{eq36}) to $g(x, z, t, u) = 1_{\{t\wedge u > s\}}$ yields
\begin{equation}\label{eq44}
n_1^- \sum_{i \in G^1} Y^i_{s+} \to^p \widetilde{q}_s \coloneqq E[Y^1_{s+} \mid X^1 \in \calx, Z^1 = 1].
\end{equation}
The function $\widetilde{q}$ is continuous on $[0, \infty)$ and 
it holds that $\widetilde{q}_s = E[Y^1_s \mid X^1 \in \calx, Z^1 = 1]$ for all $s \in [0, \tau]$.
Indeed, by the dominated convergence theorem, we have
\begin{align*}
\widetilde{q}_{s+}
&= \lim_{h \downarrow 0} E[Y^1_{(s+h)+} \mid X^1 \in \calx, Z^1 = 1] \\
&= \lim_{h \downarrow 0} E[1_{\{\widetilde{T}^1 > s+h \}} \mid X^1 \in \calx, Z^1 = 1] \\
&= E[1_{\{\widetilde{T}^1 > s \}} \mid X^1 \in \calx, Z^1 = 1] \\
&= E[Y^1_{s+} \mid X^1 \in \calx, Z^1 = 1] \\
&= \widetilde{q}_{s}
\end{align*}
for all $s \in [0, \tau]$.
Thus, $\widetilde{q}$ is right-continuous on $[0, \infty)$.
By the same way, for all $s \in (0, \tau]$, we have
\begin{align*}
\widetilde{q}_{s-}
&= \lim_{h \downarrow 0} E[Y^1_{(s-h)+} \mid X^1 \in \calx, Z^1 = 1] \\
&= \lim_{h \downarrow 0} E[1_{\{\widetilde{T}^1 > s-h \}} \mid X^1 \in \calx, Z^1 = 1] \\
&= E[1_{\{\widetilde{T}^1 \ge s \}} \mid X^1 \in \calx, Z^1 = 1] \\
&= E[Y^1_s \mid X^1 \in \calx, Z^1 = 1].
\end{align*}
Therefore, for all $s \in [0, \tau]$, it follows from $[B4]$ that
\[
\Delta \widetilde{q}_s = - E[1_{\{\widetilde{T}^1 = s \}} \mid X^1 \in \calx, Z^1 = 1] = 0
\]
and hence
\[
\widetilde{q}_s = \widetilde{q}_{s-} = E[Y^1_s \mid X^1 \in \calx, Z^1 = 1].
\]

Since $n_1^-\sum_{i \in G^1} Y^i_+$ is decreasing and $\widetilde{q}$ is continuous, 
it follows from (\ref{eq44}) and Theorem 3.37 (a) in Chapter 6 of \cite{jacod2003limit} that 
$n_1^-\sum_{i \in G^1} Y^i_+ \to^d \widetilde{q}$ in $D([0, \tau])$.
In particular, $n_1^-\sum_{i \in G^1} Y^i_+$ is C-tight and hence (\ref{eq47}) holds.
The proof of the (C-)tightness of $(\bbV^n)$ is complete.

For all $s \in [0, \tau]$, 
by ($\ref{eq36}$) with $g(x, z, t, u) = 1_{\{t\wedge u > s\}} p(s \mid x, z)$, 
we have
\[
\bbV^n_s \to^p q_s = E[Y^1_{s+}p(s \mid X^1, Z^1) \mid X^1 \in \calx, Z^1 = 1].
\]
We can show that $q$ is continuous on $[0, \tau]$ and $q_s = E[Y^1_sp(s \mid X^1, Z^1) \mid X^1 \in \calx, Z^1 = 1]$ 
in the same way as for $\widetilde{q}$.
Since $(\bbV^n - q)$ is tight, $\bbV^n - q \to^d 0$ in $D([0, \tau])$ holds.
By the continuous mapping theorem, we obtain 
\[
\sup_{s \in [0, \tau]} |\bbV^n_s - q_s| \to^p 0.
\]
From this and 
\[
\left| n_1^- \sum_{i \in G^1} Y^i_s p(s \mid X^i, Z^i) - q_s \right|
= |\bbV^n_{s -} - q_{s - }|
\le \sup_{t \in [0,  \tau]} |\bbV^n_t - q_t|
\]
for all $s \in (0, \tau]$, it holds that
\[
\sup_{s \in (0, \tau]} \left| n_1^- \sum_{i \in G^1} Y^i_s p(s \mid X^i, Z^i) - q_s \right|
\le \sup_{t \in [0,  \tau]} |\bbV^n_t - q_t|
\to^p 0.
\]
Moreover, we have
\[
n_1^- \sum_{i \in G^1} Y^i_0 p(0 \mid X^i, Z^i) = \bbV^n_0 \to^p q_0.
\]
Therefore, $[A2]$ holds.
\qed


\section*{Declarations}
\subsection{Funding}
This work was in part supported by 
Japan Science and Technology Agency CREST JPMJCR2115; 
the research project with Kitasato University; 
Japan Society for the Promotion of Science Grants-in-Aid for Scientific Research 
No. 
23H03354 (Scientific Research);  
and by a Cooperative Research Program of the Institute of Statistical Mathematics.

\subsection{Conflict of Interest}
The authors declare no conflicts of interest associated with this paper.

\subsection{Data Availability}
No datasets were generated or analysed during the current study.

\section{Appendix}
We prove Proposition \ref{prop1} and \ref{2209250620}.
We recall the following results, for convenience of reference. 
For a probability space $(\Omega, \calf, P)$ and sub $\sigma$-fields $\cala$, $\calb$ and $\calh$ of $\calf$, 
we write $\cala\indep_\calh\calb$ if $\cala$ and $\calb$ are $\calh$-conditionally independent, 
that is, $P(A \cap B | \calh) = P(A | \calh) P(B | \calh)$~a.s.\@ for all $A \in \cala$ and $B \in \calb$.

\begin{lemma}\label{lem1}
Let $(\Omega, \calf, P)$ be a probability space.
\bd
\item[(a)] Let $V \colon \Omega \to \mathbb{R}$ be an integrable random variable and 
let $\mathcal{G}, \mathcal{H}$ be sub $\sigma$-fields of $\mathcal{F}$.
If $\sigma[V]\indep_{\calh}\calg$, 
then $E[V | \mathcal{G}\vee \mathcal{H}] =  E[V | \mathcal{H}]$.
\im[(b)] 
Suppose that 
$\calv\indep_{\calh}\calg$ 
for sub $\sigma$-fields $\wh{\calg},\calg,\calv,\calh$ of $\calf$ such that $\wh{\calg}\subset\calg$. 
Then 
$\calv\indep_{\wh{\calg}\vee\calh}\calg$. 
\ed
\end{lemma}
\proof
(a) is obvious. 
Indeed, for all bounded $\calg$-measurable random variable $G$ and bounded $\calh$-measurable random variable $H$, 
it holds that
\[
E[VGH] = E[E[VG | \calh] H] = E[E[V | \calh] E[G \mid \calh] H] = E[E[V | \calh] GH]
\]
and we obtain (a).
 
We will show (b). Let $V,G,\wh{G},H$ be bounded random variables measurable with respect to $\calv$, $\calg$, $\wh{\calg}$, $\calh$, respectively. 
Then 
\beas 
E\big[VG\wh{G}H\big]
&=& 
E\big[E[V|\calh]E[G\wh{G}|\calh]H\big]\quad(\because\calv\indep_\calh\calg,\>\wh{\calg}\subset\calg)
\nn\\&=&
E\big[E[V|\calh]G\wh{G}H\big]
\nn\\&=&
E\big[E[V|\wh{\calg}\vee\calh]G\wh{G}H\big]\quad(\because(a))
\nn\\&=&
E\big[E[V|\wh{\calg}\vee\calh]E[G|\wh{\calg}\vee\calh]\wh{G}H\big]
\eeas
and hence $E\big[VG|\wh{\calg}\vee\calh\big]=E[V|\wh{\calg}\vee\calh]E[G|\wh{\calg}\vee\calh]$ a.s. 
This proves (b). \qed\halflineskip
%

%
%

%
\begin{lemma}\label{lem3}
We consider the same random variables as section \ref{section1} and suppose that condition $[C]$ holds.
Let $i \in \bbG$.
\bd
\item[(a)] $P^{T^i}(\cdot|\calf_0)=P^{T^i}(\cdot|X^i,Z^i)$ a.s.\@ for a regular conditional distribution.
\item[(b)] The regular conditional distribution $P^{(T^i, U^i)}(\cdot | \mathcal{F}_0)$ equals the product measure of 
$P^{T^i}(\cdot | \mathcal{F}_0)$ and $P^{U^i}(\cdot | \mathcal{F}_0)$.
\item[(c)] For any $t \ge 0$, we define sub $\sigma$-fields $\calg^{i}_t$ and $\calg^{\widehat{i}}_t$ of $\calf$ by setting
$\calg^{i}_t = \sigma[N^i_\theta, \dot{N}^i_\theta | \theta \le t]$ and 
$\calg^{\widehat{i}}_t = \sigma[N^j_\theta, \dot{N}^j_\theta | \theta \le t,\, j \not= i]$. 
Let $0 \le s \le t \le \tau$. Then, it holds that $\calg^{\widehat{i}}_t \indep_{\calg^i_s \vee \calf_0} \calg^i_t$.
\ed
\end{lemma}

\proof
Condition (a) follows from $[C]$ (iii) and Lemma \ref{lem1} (a).
Condition (b) is derived from $[C]$ (ii).
Let $\calv = \calg^{\widehat{i}}_t$, $\calg = \calg^i_t$, $\wh{\calg} = \calg^i_s$ and $\calh = \calf_0$ respectively 
in Lemma \ref{lem1} (b), we obtain (c).
\qed\halflineskip
\noindent
{\it Proof of Proposition \ref{prop1}.}
Suppose that $F\in\calf_0$, $0\leq s<t\leq\tau$ and $p_k,q_\ell\in[0,s]$ for $k\in\{1,...,K\}$ and $\ell\in\{1,...,L\}$. 
Let $\overline{p} = \max_{1 \le k \le K} p_k$ and $\overline{q} = \max_{1 \le \ell \le L} q_l$.
\bea\label{2209250650}&&
E\bigg[(N^i_t-N^i_s)\prod_{k=1}^K1_{\{N^i_{p_k}=0\}}\prod_{\ell=1}^L1_{\{\dot{N}^i_{q_\ell}=0\}}1_F\bigg]
\\&=&
E\bigg[1_{\{s<T^i\le t,\,T^i\leq U^i\}}1_{\{T^i\wedge U^i>\overline{p}\}\cup\{T^i>U^i\}}
1_{\{T^i\wedge U^i>\overline{q}\}\cup\{T^i \le U^i\}} 1_F\bigg]
\nn\\&=&
E\bigg[1_{\{s<T^i\le t,\,T^i\leq U^i\}}1_F\bigg]
\nn\\&=&
E\bigg[E\big[1_{\{s<T^i\le t,\,T^i\leq U^i\}}|\calf_0\big]1_F\bigg]
\nn\\&=&
E\bigg[\int_s^t E\big[1_{\{\theta\leq U^i\}}|\calf_0\big]P^{T^i}(d\theta|\calf_0)1_F\bigg]\quad(\because\text{Lemma \ref{lem3} (b)})
\nn\\&=&
E\bigg[\int_s^t P(U^i\geq\theta|\calf_0)P^{T^i}(d\theta|\calf_0)1_F\bigg]. 
\nn
\eea
On the other hand, 
\bea\label{2209250651} &&
E\bigg[(A^i_t - A^i_s)
\prod_{k=1}^K1_{\{N_{p_k}=0\}}\prod_{\ell=1}^L1_{\{\dot{N}_{q_\ell}=0\}}1_F\bigg]
\\&=&
E\bigg[\int_s^t1_{\{T^i\wedge U^i\geq\theta\}}\frac{P^{T^i}(d\theta|\calf_0)}{P^{T^i}([\theta, \infty)|\calf_0)}
1_{\{T^i\wedge U^i>\overline{p}\}\cup\{T^i>U^i\}}
1_{\{T^i\wedge U^i>\overline{q}\}\cup\{T^i \le U^i\}} 1_F\bigg]
\nn\\&=&
E\bigg[\int_s^t1_{\{T^i\wedge U^i\geq\theta\}}\frac{P^{T^i}(d\theta|\calf_0)}{P^{T^i}([\theta, \infty)|\calf_0)}1_F\bigg]
\nn\\&=&
E\bigg[\int_{\bbR_+^2}
\int_s^t1_{\{\alpha\wedge\beta\geq\theta\}}\frac{P^{T^i}(d\theta|\calf_0)}{P^{T^i}([\theta, \infty)|\calf_0)}\,
P^{(T^i, U^i)}(d\alpha d\beta|\calf_0) 1_F\bigg]
\nn\\&=&
E\bigg[\int_s^tP(T^i\wedge U^i\geq\theta|\calf_0)\frac{P^{T^i}(d\theta|\calf_0)}{P^{T^i}([\theta, \infty)|\calf_0)}1_F\bigg]
\quad(\because\text{Fubini's theorem})
\nn\\&=&
E\bigg[\int_s^tP(T^i\geq\theta|\calf_0)P(U^i\geq\theta|\calf_0)\frac{P^{T^i}(d\theta|\calf_0)}{P^{T^i}([\theta, \infty)|\calf_0)}1_F\bigg]\quad(\because\text{Lemma \ref{lem3} (b)})
\nn\\&=&
E\bigg[\int_s^tP(U^i\geq\theta|\calf_0)P^{T^i}(d\theta|\calf_0)1_F\bigg].
\nn
\eea
Therefore, the two expressions (\ref{2209250650}) and (\ref{2209250651}) coincide.
If $K=L=0$ and $F=\Omega$, we obtain $E[A^i_t-A^i_s] = E[N^i_t-N^i_s] < \infty$.
As a matter of fact, more generally, 
\beas 
E\bigg[(N^i_t-N^i_s)\prod_{k=1}^K1_{\{N^i_{p_k}=\ep_k\}}\prod_{\ell=1}^L1_{\{\dot{N}^i_{q_\ell}=\eta_\ell\}}1_F\bigg]
&=&
E\bigg[(A^i_t-A^i_s)\prod_{k=1}^K1_{\{N_{p_k}=\ep_k\}}\prod_{\ell=1}^L1_{\{\dot{N}_{q_\ell}=\eta_\ell\}}1_F\bigg]
\eeas
for any $\ep_1,...,\ep_K,\eta_1,...,\eta_L\in\{0,1\}$. 
Therefore, 
\beas 
E\big[N^i_t-N^i_s \big| \calg^i_s \vee \calf_0\big] 
\yeq
E\big[A^i_t-A^i_s \big| \calg^i_s \vee \calf_0\big]
\quad a.s.
\eeas
Consequently, Lemma \ref{lem3} (c) and Lemma \ref{lem1} (a) implies
\beas 
E[N^i_t-N^i_s|\calf_s] 
\yeq
E\big[A^i_t-A^i_s|\calf_s\big]
\quad a.s.
\eeas
By definition, $(A^i_t)_{t\in[0,\tau]}$ is \cadlag $\bbF$-predictable, and hence the proof is completed. 
\qed\halflineskip

\noindent
{\it Proof of Proposition \ref{2209250620}.} 
Since $M^i$ has bounded variation and $\Delta A^i=0$, we have 
\bea\label{0409270127}
[M^i]_t = \sum_{s\leq t}(\Delta M^i_s)^2 
=  \sum_{s\leq t}(\Delta N^i_s)^2 = N^i_t
\eea
Thus, $\langle M^i \rangle = A^i$ since $[M^i]-\langle M^i \rangle$ is a local martingale. 
The local martingale $M^i=(M^i_t)_{t\in[0,\tau]}$ is square-integrable since 
$\sup_{t\in[0,\tau]}E[(M^i_t)^2]\leq E[\langle M^i \rangle_\tau] = E[A^i_\tau] < \infty$. 

Let $i, i' \in \bbG$ be distinct and $0 \le s \le t \le \tau$.
By Lemma Lemma \ref{lem3} (c) and Lemma \ref{lem1} (a), it holds that
\begin{align*}
E[M^i_t M^{i'}_t | \calf_s] 
&= E[ E[M^i_t | \calg^{\widehat{i}}_t \vee \calg^i_s \vee \calf_0] M^{i'}_t | \calf_s] \\
&= E[ E[M^i_t | \calg^i_s \vee \calf_0] M^{i'}_t | \calf_s] \\
&= E[ E[M^i_t | \calg^{\widehat{i}}_s \vee \calg^i_s \vee \calf_0] M^{i'}_t | \calf_s] \\
&= E[ E[M^i_t | \calf_s] M^{i'}_t | \calf_s] \\
&= M^i_s M^{i'}_s.
\end{align*}
This gives $\langle M^i, M^{i'} \rangle_t = 
1_{\{i=i'\}}\int_0^t(1-\Delta A^i_s)dA^i_s=1_{\{i=i'\}}A^i_t$.
%
\qed\halflineskip 


Finally, we recall Rebolledo's martingale central limit theorem.
See \cite{helland1982central} for the proof of this theorem.

\begin{theorem}\label{thm1}
For each $n \in \mathbb{N}$, let $M^n$ be a locally square integrable martingale defined on 
a filtered probability space $(\Omega^n, \mathcal{F}^n, (\mathcal{F}^n_t)_{t \in \mathbb{R}_+}, P^n)$.
Suppose that the filtration $(\mathcal{F}^n_t)_{t \in \mathbb{R}_+}$ satisfies the usual condition.
Let $\alpha \colon \mathbb{R}_+ \to \mathbb{R}$ be a continuous function and 
$(\mathbb{W}_t)$ be a process with independent increments, $E[\mathbb{W}_t] = 0$ and
$E[\mathbb{W}_t \mathbb{W}_s] = \alpha(s \wedge t)$.
For $\epsilon > 0$, we define the adapted increasing process $\sigma^\epsilon [M^n]$ by setting 
\[
\sigma^\epsilon [M^n]_t = \sum_{s \le t} |\Delta M^n_s|^2 1_{\{ |\Delta M^n_s| > \epsilon \}} \qquad (t \in \bbR_+).
\]
Note that $\sigma^\epsilon [M^n]$ is locally integrable 
since $\sigma^\epsilon [M^n]_t \le \sigma^\epsilon [M^n]_{t-} + |\Delta M^n_t|^2$ and 
$M^n$ is a locally square integrable martingale.
In particular, we can consider the compensator of the process $\sigma^\epsilon [M^n]$.
Suppose that the following conditions hold.
\begin{enumerate}[label = (\roman*)]
\item $M^n_0 \to^p 0$.
\item $\langle M^n \rangle_t \to^p \alpha(t)$ for all $t \in [0, \tau]$.
\item $\sigma^\epsilon[M^n]^p_\tau \to^p 0$.
\end{enumerate}
Then, we have
\[
M^n \to^d \mathbb{W} \text{ on } D[0, \tau]
\] where $D[0, \tau]$ is the space of all c\`{a}dl\`{a}g function 
on $[0, \tau]$ equipped with Skorohod topology. 
\end{theorem}





\bibliographystyle{spbasic}
\bibliography{refs_cem+logrank} 

\begin{thebibliography}{14}
\providecommand{\natexlab}[1]{#1}
\providecommand{\url}[1]{{#1}}
\providecommand{\urlprefix}{URL }
\expandafter\ifx\csname urlstyle\endcsname\relax
  \providecommand{\doi}[1]{DOI~\discretionary{}{}{}#1}\else
  \providecommand{\doi}{DOI~\discretionary{}{}{}\begingroup
  \urlstyle{rm}\Url}\fi
\providecommand{\eprint}[2][]{\url{#2}}

\bibitem[{Abadie and Imbens(2012)}]{abadie2012martingale}
Abadie A, Imbens GW (2012) A martingale representation for matching estimators.
  Journal of the American Statistical Association 107(498):833--843

\bibitem[{Austin(2013)}]{austin2013performance}
Austin PC (2013) The performance of different propensity score methods for
  estimating marginal hazard ratios. Statistics in medicine 32(16):2837--2849

\bibitem[{Austin(2014)}]{austin2014use}
Austin PC (2014) The use of propensity score methods with survival or
  time-to-event outcomes: reporting measures of effect similar to those used in
  randomized experiments. Statistics in medicine 33(7):1242--1258

\bibitem[{Billingsley(1999)}]{billingsley1999convergence}
Billingsley P (1999) Convergence of probability measures, 2nd edn. No. .
  Probability and statistics in Wiley series in probability and mathematical
  statistics, Wiley, New York

\bibitem[{Br{\'e}maud(1981)}]{bremaud1981point}
Br{\'e}maud P (1981) Point processes and queues. Springer

\bibitem[{Chen and Tsiatis(2001)}]{chen2001causal}
Chen PY, Tsiatis AA (2001) Causal inference on the difference of the restricted
  mean lifetime between two groups. Biometrics 57(4):1030--1038

\bibitem[{Fleming and Harrington(2005)}]{fleming2005counting}
Fleming TR, Harrington DP (2005) Counting processes and survival analysis.
  Wiley-interscience paperback series, Wiley-Interscience, New Jersey

\bibitem[{Guo et~al(2021)Guo, Zhou, and Sun}]{guo2021adjusted}
Guo YM, Zhou J, Sun LQ (2021) Adjusted log-rank test with double inverse
  weighting under dependent censoring. Acta Mathematica Sinica, English Series
  37(10):1573--1585

\bibitem[{Helland(1982)}]{helland1982central}
Helland IS (1982) Central limit theorems for martingales with discrete or
  continuous time. Scandinavian Journal of Statistics pp 79--94

\bibitem[{Iacus et~al(2011)Iacus, King, and Porro}]{iacus2011multivariate}
Iacus SM, King G, Porro G (2011) Multivariate matching methods that are
  monotonic imbalance bounding. Journal of the American Statistical Association
  106(493):345--361

\bibitem[{Iacus et~al(2019)Iacus, King, and Porro}]{iacus2019theory}
Iacus SM, King G, Porro G (2019) A theory of statistical inference for matching
  methods in causal research. Political Analysis 27(1):46--68

\bibitem[{Ikeda and Watanabe(1989)}]{ikeda1989stochastic}
Ikeda N, Watanabe S (1989) Stochastic differential equations and diffusion
  processes, 2nd edn. No. v. 24 in North-Holland mathematical library,
  North-Holland Pub. Co.,Kodansha, Tokyo

\bibitem[{Jacod and Shiryaev(2003)}]{jacod2003limit}
Jacod J, Shiryaev AN (2003) Limit theorems for stochastic processes, 2nd edn.
  No. 288 in Die Grundlehren der mathematischen Wissenschaften,
  Springer-Verlag, Berlin

\bibitem[{Xie and Liu(2005)}]{xie2005adjusted}
Xie J, Liu C (2005) Adjusted kaplan--meier estimator and log-rank test with
  inverse probability of treatment weighting for survival data. Statistics in
  medicine 24(20):3089--3110

\end{thebibliography}

\end{document}